\documentclass[11pt]{article}

\usepackage{amsmath, amsfonts,amssymb, amsthm} \usepackage{euscript,makeidx,color,mathrsfs,latexsym}
\usepackage{graphicx,subfigure}
\usepackage{url}
\usepackage{booktabs}
\usepackage{multirow}

\newtheorem{remark}{Remark}

\newtheorem{assumption}{Assumption}

\newcommand{\EE}{\mathbb{E}}
\newcommand{\PP}{\mathbb{P}}

\newcommand{\ba}{\begin{array}}
\newcommand{\ea}{\end{array}}
\newcommand{\be}{\begin{equation}}
\newcommand{\ee}{\end{equation}}
\newcommand{\bee}{\begin{equation*}}
\newcommand{\eee}{\end{equation*}}

\DeclareMathOperator*{\argmax}{argmax}

\def\dbF{\mathbb{F}}

\def\dbP{\mathbb{P}}

\def\cF{{\cal F}}

\def\cL{{\cal L}}

\def\cT{{\cal T}}
\def\cU{{\cal U}}

\def\cX{{\cal X}}

\def\no{\noindent}

\def\bu{{\bf u}}
\def\bx{{\bf x}}

\newcommand{\basa}{\begin{assumption}}
\newcommand{\easa}{\end{assumption}}

\newcommand{\bas}{\begin{assum}}
\newcommand{\eas}{\end{assum}}

\def\bx{{\bf x}}

\def\1{{\bf 1}}

\def\:{\!:\!}

\def \proof{{\noindent \bf Proof\quad}}

at 9pt

\begin{document}

\newtheorem{thm}{Theorem}[section]
\newtheorem{lem}[thm]{Lemma}
\newtheorem{cor}[thm]{Corollary}
\newtheorem{prop}[thm]{Proposition}
\newtheorem{rem}[thm]{Remark}
\newtheorem{eg}[thm]{Example}
\newtheorem{defn}[thm]{Definition}
\newtheorem{assum}[thm]{Assumption}

\renewcommand {\theequation}{\arabic{section}.\arabic{equation}}
\def\thesection{\arabic{section}}

\title{\bf DeepPAAC: A New Deep Galerkin Method for Principal-Agent Problems}

\author{Michael Ludkovski\thanks{Department of Statistics and Applied Probability, University of California, Santa Barbara. Building 479, UC Santa Barbara, Santa Barbara, CA 93106-3110, USA. Email: ludkovski@pstat.ucsb.edu.} \qquad
Changgen Xie\thanks{Fintech Thrust, Hong Kong University of
Science and Technology (Guangzhou), Guangzhou, Guangdong Province 511453, China.  Email: cxie073@connect.hkust-gz.edu.cn.}
\qquad
Zimu Zhu\thanks{ Fintech Thrust, Hong Kong University of
Science and Technology (Guangzhou), Guangzhou, Guangdong Province 511453, China. Email: zimuzhu@hkust-gz.edu.cn.}
}
\date{\today}
\maketitle

\begin{abstract}
We consider numerical resolution of principal-agent (PA) problems in continuous time.  We formulate a generic PA model with continuous and lump payments and a multi-dimensional strategy of the agent. To tackle  the resulting Hamilton-Jacobi-Bellman equation with an implicit Hamiltonian we develop a novel deep learning method: the Deep Principal-Agent Actor Critic (DeepPAAC) Actor-Critic algorithm. DeepPAAC is able to handle multi-dimensional states and controls, as well as constraints. We investigate the role of the neural network architecture, training designs, loss functions, etc.~on the convergence of the solver,  presenting five different case studies.

\end{abstract}


\no{\it 2020 AMS Mathematics subject classification: 91B43, 68T07, 93E20. }

\section{Introduction}
Despite the importance of the Principal-Agent (PA) problem, there is a notable lack of relevant numerical algorithms. In the context of the PA problem, a significant subset of problems not only lack analytical solutions but feature challenging nonlinear partial differential equations. This necessitates the development of a general numerical algorithm to address the PA problem. 
In this article, we investigate numerical methods for continuous-time PA problems using the deep learning approach. We develop an actor-critic algorithm that can handle multi-dimensional principal-agent problems, nesting existing setups. Compared to existing setups, we make several interrelated contributions. First, we develop a new variant of the actor-critic Deep Galerkin Method (DGM) that targets the PA setup and employs multiple neural networks. Our formulation is flexible and allows for handling constrained Principal/Agent controls, multidimensional states, and multidimensional controls. Second, we investigate in detail multiple practical aspects of the algorithm, providing a comprehensive ``user guide'' on implementation. To this end, we performed experiments on the architecture of the various neural networks and how to combine multiple loss functions. We also present a detailed study regarding the training design to employ for the spatio-temporal surrogate. Third, we illustrate our approach using five different case studies, which encompass existing models, as well as new versions that do not admit closed-form solutions and have hitherto never been studied. 

Since the seminal paper by Holmstr\"{o}m and Milgrom \cite{holmstrom1987aggregation}, there has been an extensive literature on continuous-time Principal-Agent (PA) problems; see for example Sch\"{a}ttler and Sung \cite{schattler1993first}, Sung \cite{sung1995linearity,sung1997corporate}, Cvitani\'{c}, Wan and Zhang \cite{cvitanic2009optimal}, Sannikov \cite{sannikov2008continuous}, Williams \cite{williams2015solvable}, etc. Typically, there are two approaches to continuous Principal-Agent problem: the Hamilton--Jacobi--Bellman  (HJB) equation approach \cite{sannikov2008continuous,Sung2025,williams2015solvable} and the  Forward-Backward Stochastic Differential Equation (FBSDE) approach \cite{cvitanic2012contract}. 

Computationally, the PA problem can be viewed as bi-level stochastic control: the Agent's problem and the Principal's problem. The Agent's problem is solved first and the resulting process(es) are then used as part of the state process(es) for the Principal's problem, as pioneered in the seminal work of Sannikov \cite{sannikov2008continuous}. In general, both the Agent and the Principal have multiple controls, leading to multi-dimensional setups. For example, following a recent paper by Sung, Zhang, and Zhu \cite{Sung2025}, we allow both the Principal and the Agent to consume; both of these are augmented to the state of the Principal's problem. In this paper, we develop an efficient and stable algorithm to numerically solve continuous-time PA models. We follow the HJB approach, formulating the Principal's value function as solutions to a certain nonlinear partial differential equation. As a benchmark, we use the PA framework developed in \cite{Sung2025} which nests multiple existing settings and serves as a good illustration for the various flavors of PA problems. 

Despite the large volume of work, strongly focused on explicit solutions, few articles investigate numerical algorithms for PA problems. Recently, several papers utilized machine learning methods to do so \cite{baldacci2019market,campbell2021deep,dayanikli2024machine}. Baldacci et al.~\cite{baldacci2019market} adopt an actor--critic reinforcement learning approach that learns the control and contract explicitly via temporal-difference updates and policy gradients; Campbell et al.~\cite{campbell2021deep} discretize an FBSDE and train a deep--BSDE solver by minimizing the terminal mismatch for the inner problem. Dayan{\i}kl{\i}  and Lauriere \cite{dayanikli2024machine} use the method of Lagrange multipliers to rewrite the cost function and  discretize the SDE under the framework of mean field games. Compared to these works, we solve the HJB equation by directly minimizing the PDE residual at space--time samples under boundary and terminal conditions, thus requiring neither time discretization nor FBSDE simulations. 

To handle the mentioned challenges of a multi-dimensional state  with a multi-dimensional control and a possibly constrained control set, we employ the deep learning approach that uses neural networks to approximate desired quantities: Deep Principal-Agent Actor Critic (DeepPAAC). To this end, we modify the Deep Galerkin Method (DGM), first developed in \cite{sirignano2018dgm}, in order to efficiently approximate multi-dimensional HJB equations with constrained controls. Introduced by Sirignano and Spiliopoulos  to solve high dimensional PDEs, DGM was extended by Al-Aradi et al.~\cite{al2019applications} to general HJB equations, for which they introduced the Policy Improvement Algorithm (PIA) inspired by actor-critic methods \cite{jacka2017policy} in machine learning. The idea of actor-critic schemes is to utilize separate neural networks for the value function and for the optimal control which allows more flexibility in designing the loss functionals and setting up the gradient descent steps. For further variants of actor-critic schemes for HJB equations, we refer to  \cite{cheridito2025deep,cohen2025neural,duan2023relaxed,duarte2024machine,zhou2021actor}, noting that this strand is highly active and new enhancements are constantly being proposed. We also mention related actor-critic solvers for stochastic games \cite{hu2024recent}.

We utilize the ideas in \cite{al2019applications,al2022extensions} and develop a new neural network architecture inspired by the swish activation function proposed in \cite{ramachandran2017searching} and a direct encoding of the terminal constraint proposed in \cite{sukumar2022exact}. Our DeepPAAC scheme solves the Principal's HJB equation by minimizing the loss function of the differential operator, first order condition of the control, as well as a penalty function of the control (for constrained Principal/Agent controls) with a mesh free approach. Numerical results in our paper show that the algorithm is stable, accurate and fast, outperforming the traditional DGM neural network. 

To summarize our contributions, DeepPAAC offers a new PIA variant that targets generic PA formulations. As demonstrated through our multiple case studies, DeepPAAC opens the door to consider many new versions of PA models, including multi-dimensional controls, constrained contracts, and settings with multiple nonlinear terms which lead to implicit Hamiltonians.  In parallel, our adaptation of deep Actor-Critic schemes offers a new testbed of independent interest to the numerical HJB community, where we provide comprehensive experiments regarding PIA training, adaptive sampling, stopping criteria and penalty function selection for constrained settings.

The rest of the paper is organized as follows: In Section \ref{sec:model}, we first review a general continuous-time Principal-Agent model introduced in \cite{Sung2025} that serves as our benchmark, in which we characterize the value function of the Principal in the PA problem as the solution of an HJB equation. Then we introduce our new DeepPAAC method to  numerically solve this HJB equation in Section \ref{sec:Imple}.
In Section \ref{sec: Multidim-case}, we apply our algorithm to multidimensional case studies and also develop an extension for solving PA problems with constraints. We conclude our paper in Section \ref{sec:Con}.  Technical proofs are collected
in the Appendix and Python code and notebooks is publicly available at \url{https://github.com/Reedcgx/DeepPAAC}.

\section{The Model}
\label{sec:model}
We follow the formulation in \cite{Sung2025} to describe a general model for continuous time Principal-Agent (PA) problems\footnote{Technical conditions are added in \cite{Sung2025} to guarantee the well-posedness of the problem.}.

Let $(\Omega,\dbF,\dbP)$ be a standard probability space supporting a standard Brownian Motion $\{B_t\}_{\{t\geq 0\}}$. Let the filtration $\{\cF_t\}_{\{0\leq t\leq T\}}$  be the augmentation of $\{B_t\}_{\{0\leq t\leq T\}}$ on a fixed time horizon $T$.  We assume the observed output process $\{X_t\}_{\{0\leq t\leq T\}}$ satisfies 
\begin{align}\label{eq:X-dyn}
dX_t=\sigma(t,X_t)dB_t=b(t,a_t,X_t)dt+ \sigma (t,X_t)dB_t^{a}, \quad X_0=x_0,\ 0\leq t\leq T,
\end{align}
where the progressively measurable process $\{a_t\}_{\{0\leq t\leq T\}}$ is the Agent's effort and $\{B_t^a\}_{\{0\leq t\leq T\}}$ is a  Brownian Motion 
under the probability measure $\dbP^{a}$: 
\begin{align*}
\Big({d\dbP^{a}\over d \dbP}\Big)=\exp \left(\int_0^T b(t,a_t,X_t)\sigma^{-1}(t,X_t)dB_t-{1\over 2}\int_0^T \big[b(t,a_t,X_t)\sigma^{-1}(t,X_t)\big]^2dt \right)
\end{align*} 
by the Girsanov Theorem.

The effort process $\{a_t\}_{\{0\leq t\leq T\}}$ of the Agent is not observable by the Principal, who only gets to see the aggregate output $\{X_t\}_{\{0\leq t\leq T\}}$. To incentivize the Agent, the Principal offers him a contract that is linked to $\{X_t\}_{\{0\leq t\leq T\}}$ and hence to $\{a_t\}_{\{0\leq t\leq T\}}$. The contract consists of two parts: a continuous payment stream $\{S_t\}_{\{0\leq t< T\}}$, and a one-time lump sum payment $\xi$ paid at terminal time $T$.  
The continuous payment stream $\{S_t\}_{\{0\leq t<T\}}$ is characterized by the two reward processes $\{\alpha_t\}_{\{0\leq t<T\}}, \{\beta_t\}_{\{0\leq t<T\}}$ and
satisfies the dynamics
\begin{align} \notag
dS_t&=\alpha_tdt+\beta_t \,dX_t\\ \label{eq:St-dynamics}
&=[\alpha_t+\beta_tb(t,a_t,X_t)]\, dt+\beta_t \sigma(t,X_t) dB_t^{a},\ S_0=s_0,\ 0\leq t\leq T.
\end{align}

The Agent's net wealth process $\{m_t^A\}_{\{0\leq t\leq T\}}$ evolves thanks to the stream of payments $dS_t$, the consumption process  $\{c_t\}_{\{0\leq t\leq T\}}$, and interest earned at rate $r_A$. The overall dynamics of $m^A$ are thus
\begin{align} \notag
\label{eq: agent-consum}
dm^A_t&= (r_Am^A_t-c_t)\, dt+dS_t\\
&=[r_Am^A_t-c_t+\alpha_t+\beta_t b(t,a_t,X_t)]\, dt+\beta_t \sigma(t,X_t) \, dB_t^{a},\ 0\leq t\leq T.
\end{align}

The Principal's wealth process $\{m_t^P\}_{\{0\leq t\leq T\}}$ consists of paying out the continuous rewards at rate $-dS_t$, consumption (also sometimes called dividends) at rate $\{D_t\}_{\{0\leq t\leq T\}}$, and earning interest at rate $r_P$, with overall  dynamics
\begin{align} \notag
dm^P_t&= (r_Pm^P_t-D_t)dt+dX_t-dS_t\\ \label{eq: principal-consum}
&= \left[r_Pm^P_t-D_t-\alpha_t+(1-\beta_t)b(t,a_t,X_t) \right]\, dt+(1-\beta_t)\sigma (t,X_t)dB_t^{a}.
\end{align}

Now we write down the two optimization problems.

\bigskip

{\bf Agent's Problem:}
 The Agent's performance functional is 
\begin{align}
J^A (t,\alpha,\beta,\xi,a,c)=\EE^{a} \left[\int_t^T U^A(s,\alpha_s,\beta_s,a_s,c_s)ds+\Phi^A(T,\xi,m^A_T) \right],\quad 0\leq t\leq T,
\end{align}
where $U^A$ and $\Phi^A$ are both increasing, concave and integrable utility functions in $\alpha,\beta,\xi$ respectively. $\EE^a$ means the expectation is under the probability $\PP^a$ mentioned above.
Given $\{\alpha_t, \beta_t\}_{\{0\leq t\leq T\}}$ and $\xi$,  the Agent maximizes his utility by choosing optimal $\{a_t\}_{\{0\leq t\leq T\}}$ and $\{c_t\}_{\{0\leq t\leq T\}}$:
\begin{align}
\label{eq:agent}
V^A(t,\alpha,\beta,\xi)=\sup_{a,c}  J^A(t,\alpha,\beta,\xi,a,c). 
\end{align}
We denote by $a^*(s,\alpha,\beta,\xi)$ and $c^*(s,\alpha,\beta,\xi)$ the optimal controls achieving the supremum above, viewed as best response functions to the Principal's choices of $(\alpha, \beta,\xi)$. By convention, if the Agent has multiple optima $a(s,\alpha,\beta,\xi)$ and $c(s,\alpha,\beta,\xi)$, it is assumed that the Agent always chooses the one that is best for the Principal.

\medskip

{\bf Principal's problem}:
The Principal's performance functional is   based on the above best reaction of the Agent and consists of similar running and terminal cost criteria:
\begin{align}
\label{eq:principal-valuefunction}
J^P(t,\alpha,\beta,\xi,D)=\EE^{a^*(s,\alpha,\beta,\xi)} \left[\int_t^T U^P(s,\alpha_s,\beta_s,D_s)ds+\Phi^P(T,\xi,X_T,m^P_T)\right],
\end{align}
where $U^P$ and $\Phi^P$ are both increasing, concave utility functions in $D,m^P$ respectively. Observe that the expectation in \eqref{eq:principal-valuefunction} is based on the feedback probability measure $\PP^{a^*(s,\alpha,\beta,\xi)}$. The Principal's problem at time $t$ is 
\begin{align}
\label{eq:principal}
\sup_{\alpha,\beta,\xi,D} J^P(t,\alpha,\beta,\xi,D),
\end{align}
subject to the Agent's participation constraint (also called individual rationality (IR)):
$$V^A(t,\alpha,\beta,\xi)\geq R_t,$$ 
for a  given threshold $R_t$.

\subsection{HJB Equations and Critic Actor Algorithm}

For the sake of simplification (by setting $\tilde{m}^A_t=e^{-r_At}m_t^A$ and $\tilde{m}^P_t=e^{-r_pt}m_t^P$), we may assume zero interest rates $r_A=r_p=0$. The PA problem constitutes a stochastic control problem; hence, its corresponding Principal’s value function admits characterization via an associated Hamilton-Jacobi-Bellman (HJB) equation. As is standard in the literature (e.g.~\cite{sannikov2008continuous}), we first solve the Agent's problem (\ref{eq:agent}) by using the comparison principle of BSDEs and stochastic maximum principle. When solving the Agent's problem, the Agent's optimal remaining utility process is characterized as a state process, which we denote as $\{W_t\}_{\{0\leq t\leq T\}}$.    In the dynamics of $\{W_t\}$, another process $\{Z_t\}_{\{0\leq t\leq T\}}$ appears from optimizing  $J^A (t,\alpha,\beta,\xi,a,c)$ over the effort process $\{a_t\}$ using the comparison principle of BSDE. An auxiliary process $\{Y_t\}_{\{0\leq t\leq T\}}$ also arises when deriving the Agent's optimal consumption process, introducing another process $\{\hat{Z}_t\}_{0\leq t\leq T}$ from the stochastic maximum principle. These two processes, $\{W_t\}$, $\{Y_t\}$, together with (\ref{eq: agent-consum}) and (\ref{eq: principal-consum}), are the four state processes needed to solve the Principal's problem. 

These states satisfy the following dynamics for $t\leq s\leq T$ \footnote{Equation \ref{eq: principal-stateprocess} is under probability measure $P^{{I_a(\alpha_r,\beta_r,{\color{blue}X_r},Y_r,Z_r)}}$, which is equivalent to equation (5.1) in \cite{Sung2025} under the original probability $P$.} :
\begin{eqnarray}
\begin{cases}
\label{eq: principal-stateprocess}
X_s=&x+\int_t^s b(r,I_a(\alpha_r,\beta_r,X_r,Y_r,Z_r),X_r)dr+\int_t^s \sigma(r,X_r)dB_r^{I_a(\alpha_r,\beta_r,X_r,Y_r,Z_r)},\cr
W_s =&w-\int_t^s U^A(r,\alpha_r,\beta_r,I_{a}(\alpha_r,\beta_r,X_r,Y_r,Z_r),I_{c}(\alpha_r,\beta_r,X_r,Y_r,Z_r))dr\\
&+\int_t^s Z_rdB_r^{I_a(\alpha_r,\beta_r,X_r,Y_r,Z_r)}, \cr
Y_s=&y+\int_t^s \hat{Z}_rdB_r^{I_a(\alpha_r,\beta_r,X_r,Y_r,Z_r)},    \cr
m^A_s=&\int_t^s [\alpha_r-I_c(\alpha_r,\beta_r,X_r,Y_r,Z_r)+\beta_r b(r,I_a(\alpha_r,\beta_r,X_r,Y_r,Z_r),X_r)]dr\cr
&+\int_t^s \beta_r\sigma(r,X_r) dB^{I_a(\alpha_r,\beta_r,X_r,Y_r, Z_r)}_r,\cr
m^P_s=&\int_t^s((1-\beta_r)b(r,I_a(\alpha_r,\beta_r,X_r,Y_r,Z_r),X_r)-D_r-\alpha_r)dr\cr
&+\int_t^s(1-\beta_r)  \sigma(r,X_r)  dB_r^{I_a(\alpha_r,\beta_r,X_r,Y_r,Z_r)},
\end{cases}
\end{eqnarray}
with terminal  constraints \footnote{ It is well noted that by equation (\ref{eq:term}), there is a certain relationship between $W_T$ and $Y_T$, which is part of the condition for the well-posedness of equation \eqref{eq: principal-stateprocess}.}
\begin{equation}
\label{eq:term}
W_T=\Phi^A (T,\xi,m^A_T),\quad {Y}_T={\partial \Phi^A\over\partial m_T^A}(T,\xi,m_T^A).
\end{equation}
Here $w$ represents the Agent's individual reservation, namely, $w=R_t$. 
We moreover re-write the Agent's optimal controls as a function of $Z$, while
$$a^*(s,\alpha,\beta,\xi)=I_a(\alpha_s,\beta_s,X_s,Y_s,Z_s), c^*(s,\alpha,\beta,\xi)=I_c(\alpha_s,\beta_s,X_s,Y_s,Z_s)$$ for some functions $I_a,I_c$ (see \cite{Sung2025} for details).

The Principal's problem starting at $t$ is now reformulated as
\begin{align*}
V^P(t,0,w,0,0)=\sup_{y}V(t,0,w,y,0,0),
\end{align*}
where $y=Y_t$ is the initial value of $\{Y_t\}$ and 
\begin{align*}
V(t,0,w,y,0,0)=&\sup_{\alpha,\beta,D,Z,\hat{Z}}\EE^{I_a(\alpha,\beta,X,Y,Z)} \Big[\int_t^T U^P(s,\alpha_s,\beta_s,D_s)ds\\
&+\Phi^P(T,\xi,X_T,m^P_T)\Big],
\end{align*}
with respect to the state processes \eqref{eq: principal-stateprocess} and terminal constraints \eqref{eq:term}. For the terminal condition on $V$, $\xi$ appearing in the term $\Phi^P(T,\xi,X_T,m^P_T)$ of the Principal's performance functional can be eliminated by solving the constraints (\ref{eq:term}), which induces $\xi=f(m_T^A, W_T, Y_T)$ for some function $f$. This generates the terminal condition of the Principal's problem of the form 
\begin{align*}
V(T,x,w,y,m^A,m^P)=\Phi^P(T,f(m^A, w, y),{\color{blue}x},m^P).
\end{align*}
In general, the terminal condition for the HJB equation \eqref{eq:principal-operator} is not straightforward and will be discussed separately in each example. By the dynamic programming principle, $V$ is the solution to the following HJB equation: 
\begin{align}
\label{eq:principal-generalHJB}
 0 & = \partial_t V(t,\bx)+\sup_\bu\left\{\cL_t^{\bu}V(t,\bx)+F(t,\bx,\bu)\right\}, \end{align}
 where $\bx=(x,w,y,m^A,m^P)$, $\bu=(\alpha,\beta,D,Z,\hat{Z})$, $F(t,\bx,\bu)=U^P(t,\alpha,\beta,D)$ and the generator is
\begin{align}
 {\cal L}^\bu_t V&=
 V_x\cdot b(t,I_a(\alpha,\beta,X,Y,Z),x) \nonumber\\ 
 &-V_w \cdot U_A \left(t,\alpha,\beta,I_a(\alpha,\beta,X,Y,Z),I_c(\alpha,\beta,X,Y,Z) \right)  \nonumber\\  
& +V_{m^A}\cdot \! [\alpha-I_{\color{blue}c}(\alpha,\beta,X,Y,Z)+\beta b(t,I_a(\alpha,\beta,X,Y,Z),x)]\nonumber\\  
&+V_{m^p}\cdot \![(1-\beta)b(t,I_a(\alpha,\beta,X,Y,Z),x)-D-\alpha]\nonumber\\ 
& +V_{xw}\sigma(t,x)Z+V_{xy}\sigma(t,x)\hat{Z}+V_{xm^A}\beta \sigma^2(t,x)+V_{xm^p}(1-\beta)\sigma^2(t,x)\nonumber\\
&+V_{wy}Z\hat{Z}+V_{wm^A} \beta Z\sigma(t,x)+V_{wm^p} Z(1-\beta)+V_{ym^{A}}\beta\hat{Z}+V_{ym^p}\hat{Z}(1-\beta)\sigma(t,x) \nonumber\\ 
&  +V_{m^Am^p}(1-\beta)\beta\sigma^2(t,x)+{1\over 2}V_{xx}\sigma^2(t,x) +{1\over2} V_{ww}Z^2\nonumber\\
&+{1\over 2} V_{yy} \hat{Z}^2+{1\over 2}V_{m^Am^A}\beta^2\sigma^2(t,x)+{1\over 2}V_{m^pm^p} (1-\beta)^2\sigma^2(t,x).
\label{eq:principal-operator}
\end{align}

\subsection{Policy Iteration}
\label{subsec:algo}

When an explicit solution for $V$ is not available, the main approach is to solve the HJB equation \eqref{eq:principal-generalHJB}. Observe that to solve for $V$, one must first determine the functional forms of $I_a,I_c$ which come from the optimal solution $a(\cdot), c(\cdot)$ to \eqref{eq:agent}. Given the nonlinearities that arise from the Hamiltonian term, this is not possible to do in closed-form. In simpler cases one can directly evaluate the supremum defining the Hamiltonian (i.e.~the Agent's problem) and hence reduce to a nonlinear PDE that can be tackled using common numerical methods, such as finite-difference or finite-element techniques. However, one of the defining features of PA problems is that the supremum over the Agent controls is in general highly non-trivial, in part due to having multiple controls to optimize over, and in part due to their complex interaction between Principal and Agent controls via $I_a,I_c$. Consequently, the local optimization over $a,c$ might not be possible to resolve, and hence we are facing a PDE with implicitly defined terms. 

In the past, such problems were viewed as very challenging, since they require \emph{multiple levels} of approximation---locally to solve the supremum and then globally to solve for $V(t,\bx)$. Below we develop a deep-learning based setup to address this aspect of PA HJB. The key idea that we rely on is to construct (overparametrized) functional approximator that replaces $V(t,\bx)$ with a surrogate $\hat{V}(t,\bx)$ that is trained through gradient descent. The approximation is global in $(t,\bx)$ (e.g.~there is no finite-differencing) and moreover carries through to the derivatives, which are evaluated using algorithmic differentiation. In turn, functional approximation allows to do the same for the optimal control $u(\cdot, \cdot)$, which is similarly approximated via a deep neural network. 

The above construction is conceptually dimension-agnostic, and moreover thanks to being overparametrized can be trained exhaustively, breaking down the overall solver into a loop consisting of doing SGD steps over training batches that are refreshed across epochs. Consequently, the method has modest memory requirements and does not rely on any temporal or spatial discretization.

To present the numerical algorithm for  \eqref{eq:principal-generalHJB}, we consider a generic multidimensional nonlinear HJB equation written as: 
\begin{equation}
\left\{
\begin{aligned}
&\partial_t V(t,\bx)+\sup_\bu\left\{\cL_t^{\bu}V(t,\bx)+F(t,\bx,\bu)\right\}=0\\
&V(T,\bx)=G(\bx),\\
\end{aligned}
\right.
\end{equation}
with unknown value function $V$ and unknown optimal feedback control $\bu \equiv u(t,\bx)$, with
$\cL^\bu V(t,\bx)\equiv \cL^{\bu(t,\bx)}V(t,\bx)$.
Here, both $\bx$ and $\bu$ are vectors in $\cX \subseteq \mathbb{R}^d$ and $\cU \subseteq \mathbb{R}^a$ respectively.

The policy improvement algorithm (PIA) is an iterative scheme that generates a sequence of approximate solutions $V_n,u_n$ for the value function and the optimal control. Starting with an initial guess $V_0, u_0$, PIA runs as a loop for $n=1,2,\ldots$:\\

1. Find a classical solution $V_n(\cdot, \cdot)$ to the PDE
\begin{equation}\label{eq:pia-v}
\left\{
\begin{aligned}
&\partial_t V(t,\bx)+\cL_t^{\bu_{n-1}}V(t,\bx)+F(t,\bx,\bu_{n-1}(t,\bx))=0\\
&V(T,{\bx})=G({\bx}).\\
\end{aligned}
\right.
\end{equation}

2. Given $V_{n}$ compute the policy improvement 
\begin{align}\label{eq:pia-u}
\bu_{n}(t,\bx)=\argmax_{\bu \in \cU} \left\{ \cL^\bu_t V_n(t,\bx)+F(t,\bx,\bu) \right\}.
\end{align}

Thus, the idea of PIA is to alternate between the ``diffusion'' (resembling the heat equation) step to update $V_{n+1}$ given $\bu_n$ and then to update $\bu_{n+1}$ given $V_{n+1}$.
In the language of machine learning, we have the actors $V_n$ and the critics $\bu_n$.  By decoupling the two steps, PIA is able to handle complex settings such as ours, where the operator ${\cal L}$ is nonlinear, where the spatial input $\bx$ is multi-dimensional, and simultaneously where the controls $\bu$ are also multi-dimensional and cannot be expressed directly in terms of $V$. 

In the PA context, the second-order diffusion terms in \eqref{eq:pia-v}-\eqref{eq:pia-u} are, in general, controlled, which makes convergence analysis very challenging. When there is no diffusion control, \eqref{eq:pia-v} is semi-linear and convergence of the PIA is quite standard, see Huang, Wang,
and Zhou \cite{huang2025convergence}, Ito, Reisinger, and Zhang \cite{ito2021neural}, Jacka and Mijatovi{\'c} \cite{jacka2017policy} , Ma, Wang, and Zhang \cite{ma2026convergence}, Kerimkulov, Siska, and Szpruch \cite{kerimkulov2020exponential,kerimkulov2021modified}, Puterman \cite{puterman1981convergence}. Extensions to the fully nonlinear case like ours are partially addressed in  
 Jacka,  Mijatovi{\'c} and Siraj \cite{jacka2017coupling},  Kerimkulov, Siska, and Szpruch \cite{kerimkulov2021modified},Wang, and Zhang \cite{ma2026convergence}, Tran, Wang, and Zhang \cite{tran2024policy}.
However, the problems we study may not satisfy the technical conditions proposed in \cite{kerimkulov2021modified,ma2026convergence,tran2024policy}, and the respective theoretical study is left for future research.

\subsection{Principal-Agent Actor Critic Algorithm}

We now propose an instance of Policy Iteration that is specially designed for Principal-Agent problems with an implicit Hamiltonian.

Denote the machine learning approximations for the value function $V$ and optimal feedback control map as $v(t,\bx,\theta^V)$ and $u(t,\bx,\theta^u)$. The $\theta$'s are the  hyperparameters which are to be optimized as part of surrogate construction.  To this end we utilize the concept of policy improvement or actor-critic iterations described in PIA, to sequentially update $\theta^V, \theta^u$. We thus re-write as $V_n \equiv v(\cdot; \theta^V_n); u_n \equiv u(\cdot; \theta^u_n)$ and focus on learning/updating $\theta^V_n,\theta^u_n$ for $n=1,\ldots$. While we still utilize the logic of alternating between learning $\theta^V$ and $\theta^u$ as in \eqref{eq:pia-v}-\eqref{eq:pia-u}, we no longer directly optimize any objective. Rather, we simply take a few steps based on a respective stochastic gradient descent (SGD). Thus, $\theta^V_n$ is not a minimizer of any given problem, but is an update of $\theta^V_{n-1}$ based on several iterations of SGD (i.e.~it is not run until convergence). This procedure naturally incorporates the standard learning of a ML model via gradient descent with the goal of simultaneously learning the value function and the control. Moreover, by taking just a few SGD steps per iteration we aim to dramatically speed up the algorithm, in effect making as much effort as solving a single HJB equation, while PIA by definition solves \emph{many} optimization problems. To judge the quality of a given $\theta^V_n, \theta^u_n$ approximation we rely on a global, iteration-independent loss criterion that checks whether the desired equality of the HJB equation holds (in contrast, the objective criterion in PIA, say \eqref{eq:pia-u} is by construction \emph{relative} to the given $V_n$).

As an additional adjustment,  
inspired by \cite{sukumar2022exact}, we enforce the terminal condition for the value function at every step by injecting it into the output: the network predicts a correction $v(t,\bx; \theta^V_n)$ and the model output is
\begin{align*}
    V_n(t,\bx) = G(\bx) + (T - t)^{\gamma}\cdot v(t,\bx; \theta^V_n).
\end{align*}
(In our implementation we set $\gamma \equiv 1$). This reduces potential conflicts between the interior and terminal objectives in multi-objective training.

Our Principal-Agent Actor Critic (PAAC) Algorithm is as follows:

1.  Initialize value function hyperparameters $\theta_0^V$ and those for the feedback control $\theta_0^u$. Pick the learning rate schedule $(\alpha_n)$.

2. Generate training samples  ($t_n^m,\bx_n^m$), $m=1,\ldots,M$ from the domain's interior $(0,T)\times \cX$ according to the training distribution $\nu(\cdot)$.

3. Calculate the loss function (note that the sums over training sites approximate integrals over the training domains) as
    $$L_{n}\left(\theta_n^V\right)=\sum_{m=1}^{M} L_{int}\left(\theta_n^V;t_n^m,\bx_n^m \right)$$
   where the pointwise loss functional for the current epoch (based on the full $V$) is
$$L_{int}(\theta_n^V;t_n^m,\bx_n^m)=\left[\left(\partial_t+\cL^{\bu(t_n^m,\bx_n^m;\theta_n^{u})} \right){\color{blue}V}(t_n^m,\bx_n^m;\theta_n^V)+F(t_n^m,\bx_n^m,\bu(t^m_n,\bx^m_n;\theta_n^u))\right]^2.$$

4. Take a descent step for $\theta^V$  with learning rate $\alpha_n$:
\begin{align}\label{eq:sgd}
\theta_{n+1}^V\simeq \theta_n^V-\alpha_n\nabla_{\theta} L_{n}(\theta_n^V).
\end{align}
Note that in the software implementation, we employ the Adam algorithm \cite{kingma2014adam} which has additional terms related to momentum, etc., beyond \eqref{eq:sgd}.

5. Calculate the control loss functional, namely the loss metric for the control NN $\bu(\cdot,\cdot; \theta^u_n)$: $$L_u(\theta_n^u)=-\sum_{m=1}^{M} \left[\cL^{\bu(t_n^m,\bx_n^m;\theta_n^u)}{\color{blue}V}(t_n^m,\bx_n^m;\theta_n^V)+F(t_n^m,\bx_n^m,\bu(t_n^m,\bx_n^m;\theta_n^u)) \right] .$$

6. Take a descent step for $\theta^u$ with Adam-based learning rate $\alpha_n$:
\begin{align*}
\theta_{n+1}^u=\theta_n^u-\alpha_n \nabla_{\theta} L_u(\theta_n^u).
\end{align*}

7. Repeat steps (3)-(6) to make $B$ stochastic gradient descent updates.

8. Repeat steps (2)-(7) over new epochs (i.e.,~new training batches) until a stopping criterion is satisfied.

\section{Implementation}
\label{sec:Imple}
The idea of the PAAC scheme is to implement an actor-critic approach, alternately fitting the NN for the value function $V$ and for the control map $\bu$. The loss metric $L_{int}$ determines the quality of the fit for $V$ in the interior $t \in (0,T)$, $\bx \in \cX$, and the loss function $L_u$ the quality of $\bu(\cdot, \cdot)$. Note that $L_{int}$  should converge to zero everywhere, while $L_u$ should stabilize at a non-zero value. 
 The respective training amounts $M$ determine the accuracy of training; larger training sets will ensure smoother descent and hence faster convergence (but will of course be slower per epoch) across epochs.  The numerical experiments below were conducted on a workstation equipped with an Intel(R) Xeon(R) Platinum 8358P CPU at 2.60 GHz and an NVIDIA A40 GPU with 48 GB of memory.

\subsection{Stopping Criteria}
Due to the presence of multiple convergence criteria, determining the stopping criterion of our PAAC algorithm is nontrivial. 
We define the metric
\begin{align}\label{eq:L-ctrl}
L_{ctrl}(\theta^V, \theta^u;t,\bx):=({\partial L\over \partial u^1}, {\partial L\over \partial u^2},\cdots, {\partial L\over \partial u^a})^{u(t,\bx;\theta^{u})} v(t,\bx;\theta^V)
\end{align}
to verify whether step 6 in our PAAC algorithm is converging (note that if $\bu$ is vector, then $L_{ctrl}(\theta^V, \theta^u;t,\bx)$ is also a vector).
To define a stopping criterion, we construct a validation set $\{ \bar{t}^m, \bar{\bx}^m \}, m=1,\ldots, M_V$. Unlike the training batches that vary across steps, the validation set is fixed throughout, to allow an objective tracking of surrogate quality and algorithm convergence and to mitigate over-fitting to a particular training epoch. For the stopping criterion, we seek to control the $L_\infty$ norm of the vectors $\bar{L}_{int}$ and $\bar{L}_{ctrl}$ which are counterparts of $L_{int}, {L}_{ctrl}$ from \eqref{eq:L-ctrl} over the validation set. Thus, we continue the PAAC iterations until we ensure that the maximum pointwise errors of $v_n,u_n$  in the interior  are small: 
$$
\| \bar{L}_{int}( \theta^V_n; {\bar{t}}^\cdot, {\bar{\bx}}^\cdot) \|_{\infty}\le \mathrm{Tol_{int}}, \qquad \| \bar{L}_{ctrl}( \theta^V_n, \theta^u_n; {\bar{t}}^\cdot, {\bar{\bx}}^\cdot) \|_{\infty} \le \mathrm{Tol_{ctrl}}, 
$$ 
for given thresholds $\mathrm{Tol_{int}},\mathrm{Tol_{ctrl}}$. In the examples below $\mathrm{Tol_{int}}=10^{-2}$ or $10^{-3}$, $\mathrm{Tol_{ctrl}}=10^{-2}$ or $10^{-3}$. 
    
\subsection{Neural Network Architecture}{\label{sec:NNA}}

 We consider two different architectures for the neural networks. The first one is the DGM (Deep Galerkin Method) architecture used by Sirignano and Spiliopoulos \cite{sirignano2018dgm}. 
 The second one is a residual fully-connected feed-forward neural network. 
 
\begin{figure}[htbp]
\centering
\label{fig:arch.}
\subfigure[NN for value function]{\includegraphics[width=10cm]{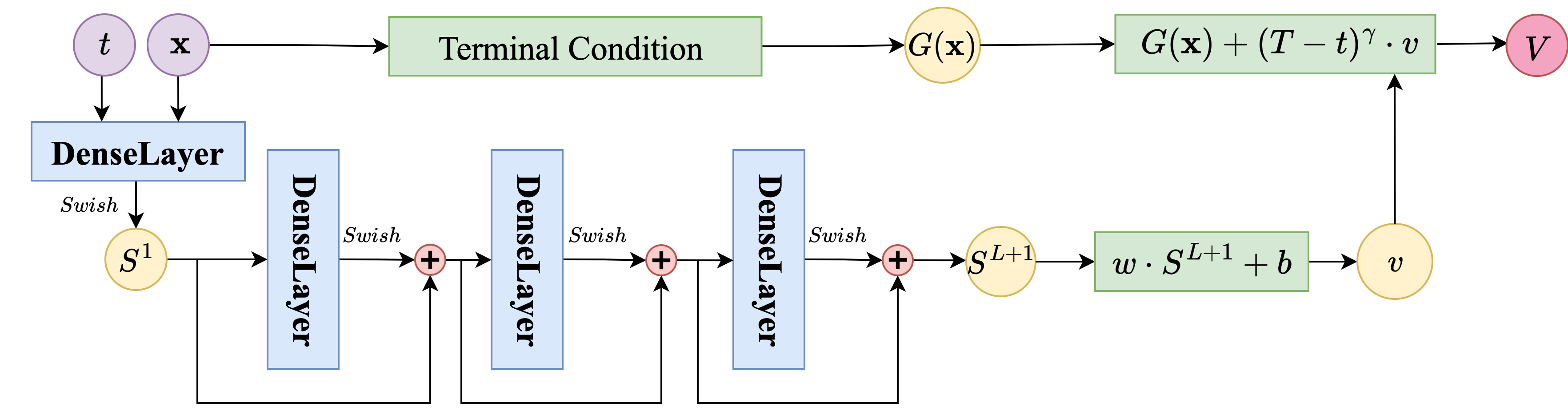}} \\
\subfigure[NN for control]{\includegraphics[width=10cm]{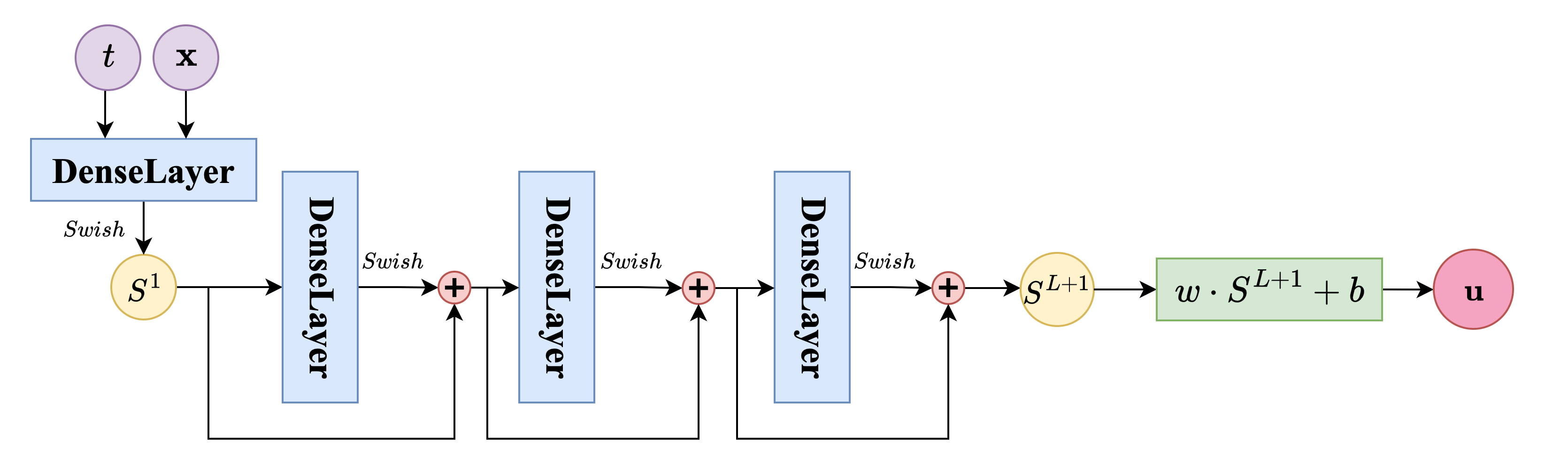}}
\caption{Top: NN architecture for value function. Below: NN architecture for control}
\end{figure}

 Unlike the DGM architecture of \cite{sirignano2018dgm}, we remove the complex LSTM-like gating and replace it with a stack of dense layers with residual connections. Motivated by \cite{ramachandran2017searching}, we use the smoother \textit{swish} activation function to improve optimization stability:
 $$Swish(x):= {x\over 1+e^{-x}}.$$
 The detailed architecture of the neural network is summarized in Figure \ref{fig:arch.}.

 To evaluate the differential operator ${\cal L}^{\bu(\cdot; \theta^u_n)}$ that defines the loss $L_{int}$ we compute numerical gradients of the current value function approximation $v(\cdot; \theta^V_n)$ using backpropagation. To this end we rely on the native backpropagation (aka algorithmic differentiation) support in modern ML libraries, such as  the function \texttt{tf.gradients} in Python's \texttt{Tensorflow} library. This capability together with the underlying universal approximation property of deep NN's are the motivating factors for using NN approximators, as it allows to replace ``term-for-term" the HJB equation with a corresponding neural network approximation.

Another issue to consider is how to handle multi-dimensional controls. Al Aradi et al.~\cite{al2019applications} suggested to use a single NN with multiple outputs, meaning that all the components of $\bu$ share the same hidden layers except for the last one, where one sets \texttt{nOutput=a}. This links all the control surrogates, aiming to exploit similarities in how the different control components depend on the system state. Alternatively, one can set up $a$ separate independently trained NN's for each control. In either situation, one has to take separate gradient descent steps for each control coordinate. This gradient descent can be done in parallel across $a$, or one-at-a-time, i.e.~sequentially.

 \subsection{Training Designs}
The training points of each batch determine what the NNs are learning. Therefore, their selection via the sampling density $\nu$ is an essential part of achieving an efficient implementation. Constructing $\bx^n_t$ becomes more complicated in multiple dimensions, where geometric intuition is trickier. One standard choice is to utilize \emph{space-filling} batch, which can be straightforwardly implemented by sampling each coordinate of $\bx$ as i.i.d. Uniform within a given range.

Non-uniform training can be used to preferentially focus on some regions.  For instance, NNs tend to have worse fit at the edges of the domain, compared to its center, so one may put more training points for $x_k$ very small or very large. Similarly, learning $V$ for $t$ small is generally harder, since the terminal condition automatically takes over for $t \to T$. Consequently, it might be beneficial to have a training density that decreases in $t$.

As a third option, to enhance convergence, we consider adaptive sampling, where new training points are picked based on the current loss $L_{int}(\cdot, \cdot)$, targeting regions where the current loss is highest. This achieves two purposes. First, the use of additional training locations helps to achieve global minimization of the loss metric. Recall that $L_n$ is supposed to represent the integrated error over the entire domain; the approximation using a fixed set of $(t^m,\bx^m)_{m=1}^M$ runs the danger of overfitting at those specific sites. Second, by preferentially picking $(t^m,\bx^m)$ we can speed-up convergence by focusing on regions where the present fit is less good.

\subsection{A toy example: The Holmstr\"{o}m--Milgrom Model}
\label{subsec:HM}
As a first illustration of our proposed numerical PAAC algorithm we consider a toy case where an exact solution is known. Namely, we take a special case of the setting in Section \ref{sec:model}: the PA model of Holmstr\"{o}m and Milgrom \cite{holmstrom1987aggregation}. Thus, we consider a contract consisting  only of a lump-sum terminal payment and assume exponential utility for the Agent. In the notation of Section \ref{sec:model},  we take  $b(t,a,x)=a, \sigma(t,x)=1$, zero out the consumption processes for the Principal and the Agent, i.e., $c_t \equiv D_t\equiv 0$, and moreover set $m^A_t= m^P_t \equiv 0$ and $U^A(\cdot)=U^P(\cdot)=0$. To be precise:

The {\bf Agent's problem} is
\begin{align}
V^A(t,\xi)=\sup_{a} \EE^a \left[\Phi^A(\xi-{1\over 2}\int_t^T a_s^2ds) \mid {\cal{F}}_{t} \right]  
\end{align}
with exponential utility $\Phi^A(x)=-\exp(-\gamma_Ax)$ for some positive $\gamma_A > 0$. We denote by $a^*(s,\xi)$ the best response in the Agent's problem. Solving the latter (see proof of Proposition \ref{prop:HMmodel}), we have the terminal payment expressed as 
\begin{align}
\label{eq:HM-Agent}
     \xi&=w+\int_t^T {\gamma_A-1\over 2}Z_s^2 \, ds+\int_t^T Z_s dB_s,
\end{align}
where $a(s,\xi)=Z_s, \, s \in [t, T]$ and $w$ is a scalar to be  determined in the Principal's problem.

The corresponding {\bf Principal's Problem} with a utility function $\Phi^P(x)$ is (without loss of generality, we assume $X_0=0$ and set $x=-w$)
\begin{align}
V^P(t,x)&=\sup_{\xi} \EE^{a(s,\xi)} \left[\Phi^P(X_T^a-\xi) \mid {\cal{F}}_t \right]\nonumber\\
&=\sup_{Z_s} \EE^{Z}\left[\Phi^P \Bigl(-w+\int_t^T (Z_s-{1+\gamma_A\over 2}Z_s^2) ds+\int_t^T (1-Z_s)dB_s^Z \Bigr) \mid {\cal{F}}_t \right]\nonumber\\
&=\sup_{Z_s} \EE^{Z}\left[\Phi^P \Bigl(x+\int_t^T (Z_s-{1+\gamma_A\over 2}Z_s^2) ds+\int_t^T (1-Z_s)dB_s^Z \Bigr) \mid {\cal{F}}_t \right], 
\end{align}
subject to (\ref{eq:HM-Agent}) and $V^A(t,\xi)\geq R_t$. When Agent's individual reservation is $R_t$, we have $w=(\Phi^A)^{-1}(R_t)=-{1\over \gamma_A}\log(-R_t)$. The Agent's optimal process (\ref{eq:HM-Agent}) is used as a state process (therefore, the Agent's HJB equation disappears) for the Principal's problem (as in Sannikov \cite{sannikov2008continuous}). We may correspondingly define the feedback control $Z(t,x)$ (similarly, $a(t,x)$). We fully solve the Principal's problem in the following proposition.

\begin{prop}
\label{prop:HMmodel}
The Principal's value function $V^P(t,x)$ can be characterized as the solution of the following HJB equation:
\begin{equation}
\label{eq:HMHJB}
\left\{
\begin{aligned}
&V^P_t+\sup_Z\left\{V^P_x(Z-{1+\gamma_A\over 2}Z^2)+{1\over 2}V^P_{xx}(1-Z)^2\right\} =0\\
&V^P(T,x) = \Phi^P(x).\\
\end{aligned}
\right.
\end{equation}
When $\Phi^P(x)=-\exp(-\gamma_P x), \gamma_P >0$ is of exponential type, we have the explicit solution:
\begin{align}
\label{eq: HMsol}
V^P(t,x)=\Phi^P \left(x+\Bigl[{1\over 2}{(1+\gamma_P)^2\over 1+\gamma_A+\gamma_P}-{\gamma_P\over 2}\Bigr](T-t) \right),
\end{align}
and Agent's optimal control is constant
\begin{align}\label{eq:a-HM}
a^*(t,x)=Z^*(t,x)={V^P_x-V^P_{xx}\over V^P_x(1+\gamma_A)-V^P_{xx}} \equiv{ 1+\gamma_P \over 1+\gamma_A+\gamma_P}.
\end{align}
The optimal terminal lump-sum contract is 
\begin{align}\label{eq:xi-HM}
\xi^*=-{1\over\gamma_A}\log(-R_t)+{\gamma_A-1\over 2}\left({1+\gamma_P\over 1+\gamma_A+\gamma_P} \right)^2(T-t)+{1+\gamma_P\over 1+\gamma_A+\gamma_P}(B_T-B_t).
\end{align}
\end{prop}
\proof See Appendix A. 

\begin{remark}
The static problem (from time $0$ to time $T$) of Holmstr\"{o}m and Milgrom's model has been well studied, see \cite{cvitanic2009optimal,holmstrom1987aggregation,schattler1993first}, etc. However, the HJB equation (\ref{eq:HMHJB}) and the solution (\ref{eq: HMsol}) for the dynamic problem, to the best of our knowledge, is new for the literature.  
\end{remark}

For this example,  taking the gradient of the Hamiltonian we have the convergence loss criterion 
\begin{align} \label{eq:foc-3.4}
L_{ctrl}=\Big[V_x^P[1-(1+\gamma_A)Z]+V_{xx}^P(Z-1)\Big].
\end{align}
For a numerical illustration we choose the mixture-of-exponentials utility
\begin{align}
\Phi^P(x)=-\lambda e^{-\gamma_P x}-(1-\lambda)e^{-\tilde{\gamma}_Px}
\end{align}
for some positive numbers $\gamma_P$,$\tilde{\gamma}_P$ and mixture weight $\lambda\in[0,1]$.
In this case, the value function and the optimal control for the HJB equation (\ref{eq:HMHJB}) have no explicit solution. However, from Proposition \ref{prop:HMmodel}, we know that if $\Phi^P(x)=-e^{-\gamma_P x}$, then the HJB equation has an explicit solution \eqref{eq: HMsol}, so we can compare to the respective linear interpolation approximation: 
\begin{align}
a^{lin}(t,x):=  \lambda {1+\gamma_P\over 1+\gamma_A+\gamma_P}+(1-\lambda){1+\tilde{\gamma}_P\over 1+\gamma_A+\tilde{\gamma}_P}.
\end{align}
Note that when $\lambda \in \{0,1\}$ the optimal strategy is a constant, independent of $(t,x)$, and so is $a^{lin}$. This is not the case for $\lambda \in (0,1)$.

\begin{figure}[htbp]
\begin{center}
\includegraphics[height=2.2in,trim=0.4in 0.6in 0.4in 0.3in]{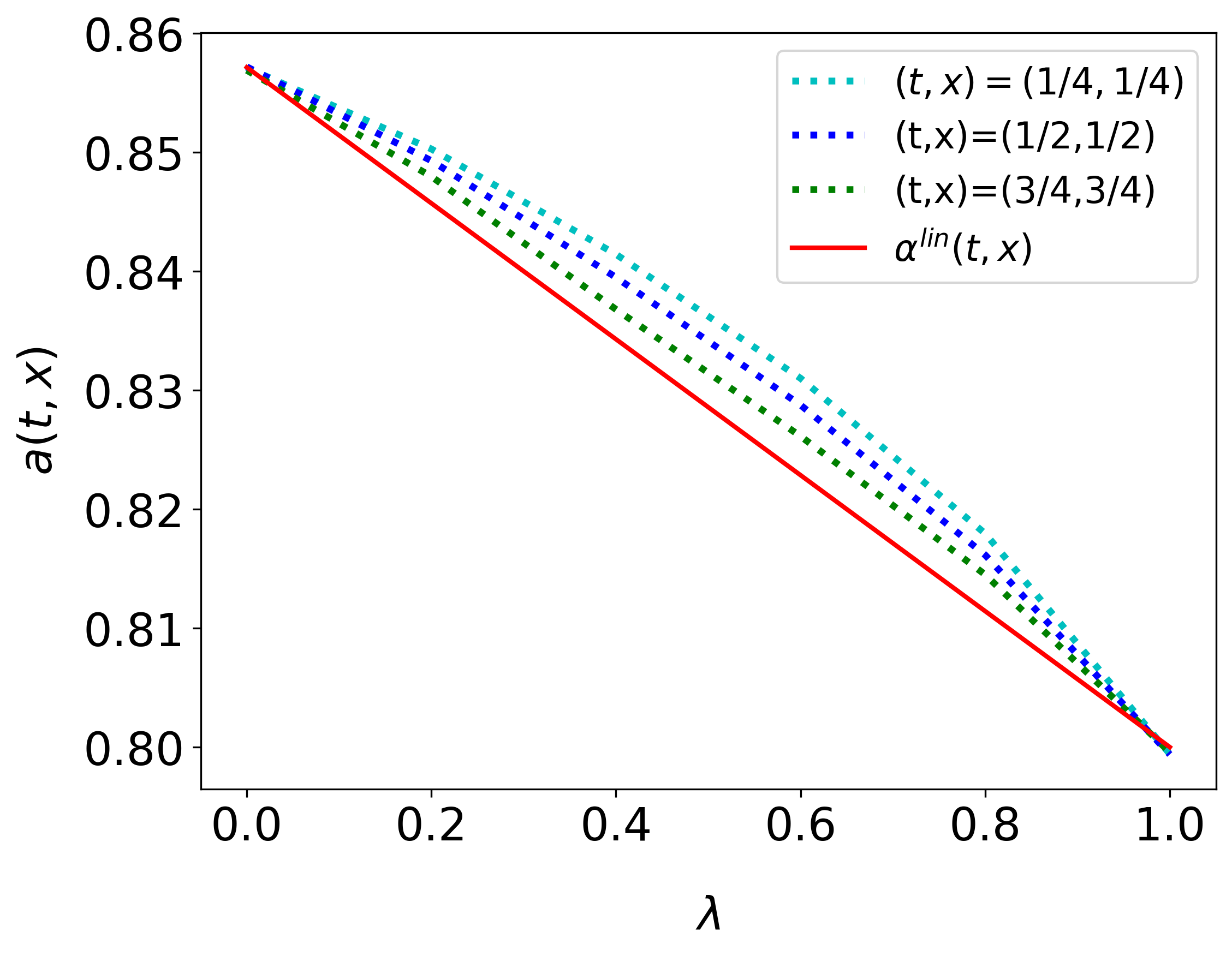}
\end{center}
\caption{ Estimated optimal control ${a}(t,x)$ obtained from the DeepPAAC algorithm for different exponential mixture weights $\lambda \in (0,1)$.}\label{fig:control for different c}
\end{figure}

\begin{figure}[htbp]
\centering
\subfigure[Relative error of ${V^P(t,x;0)}$ for $\lambda=0$]{\includegraphics[width=5cm,trim=1.4cm 0.1cm 1.4cm 2.0cm]{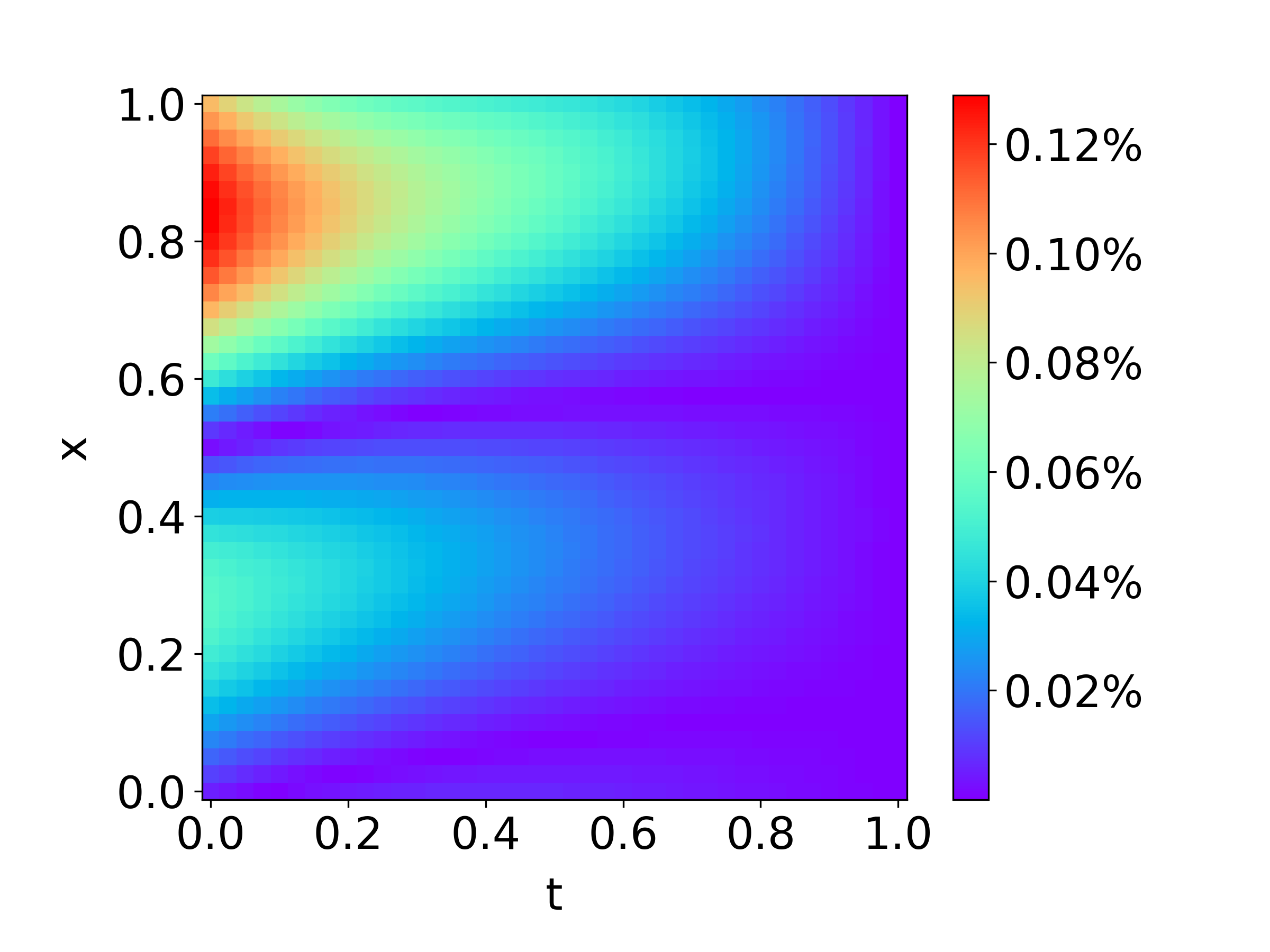}} 
\subfigure[${a}(t,x;0)$ for $\lambda=0$]{\includegraphics[width=5cm,trim=0.1cm 0.1cm 0.1cm 0.2cm]{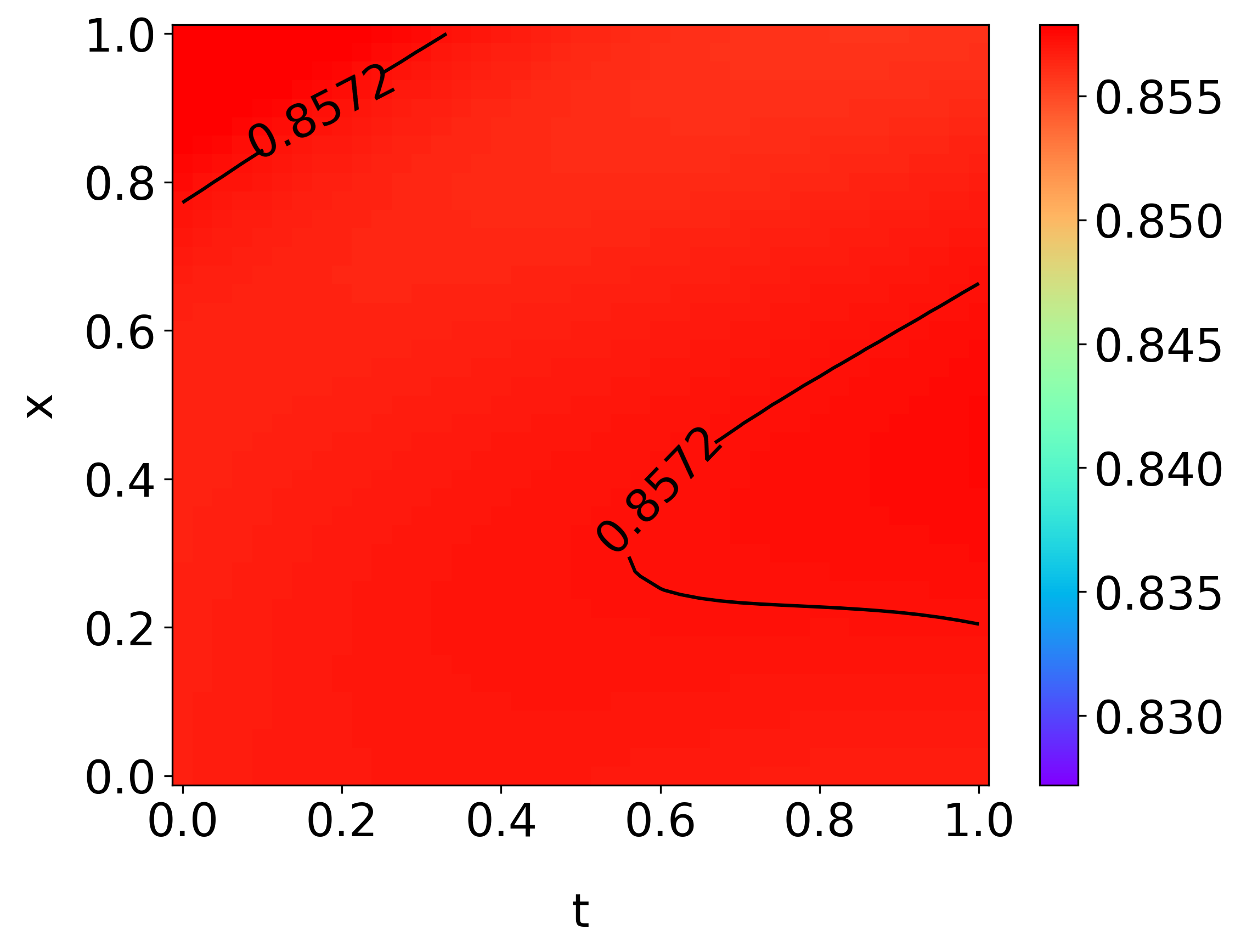}} 
\subfigure[${a}(t,x;0)$ for $\lambda=0.5$]{\includegraphics[width=5cm,trim=0.1cm 0.1cm 0.1cm 0.2cm]{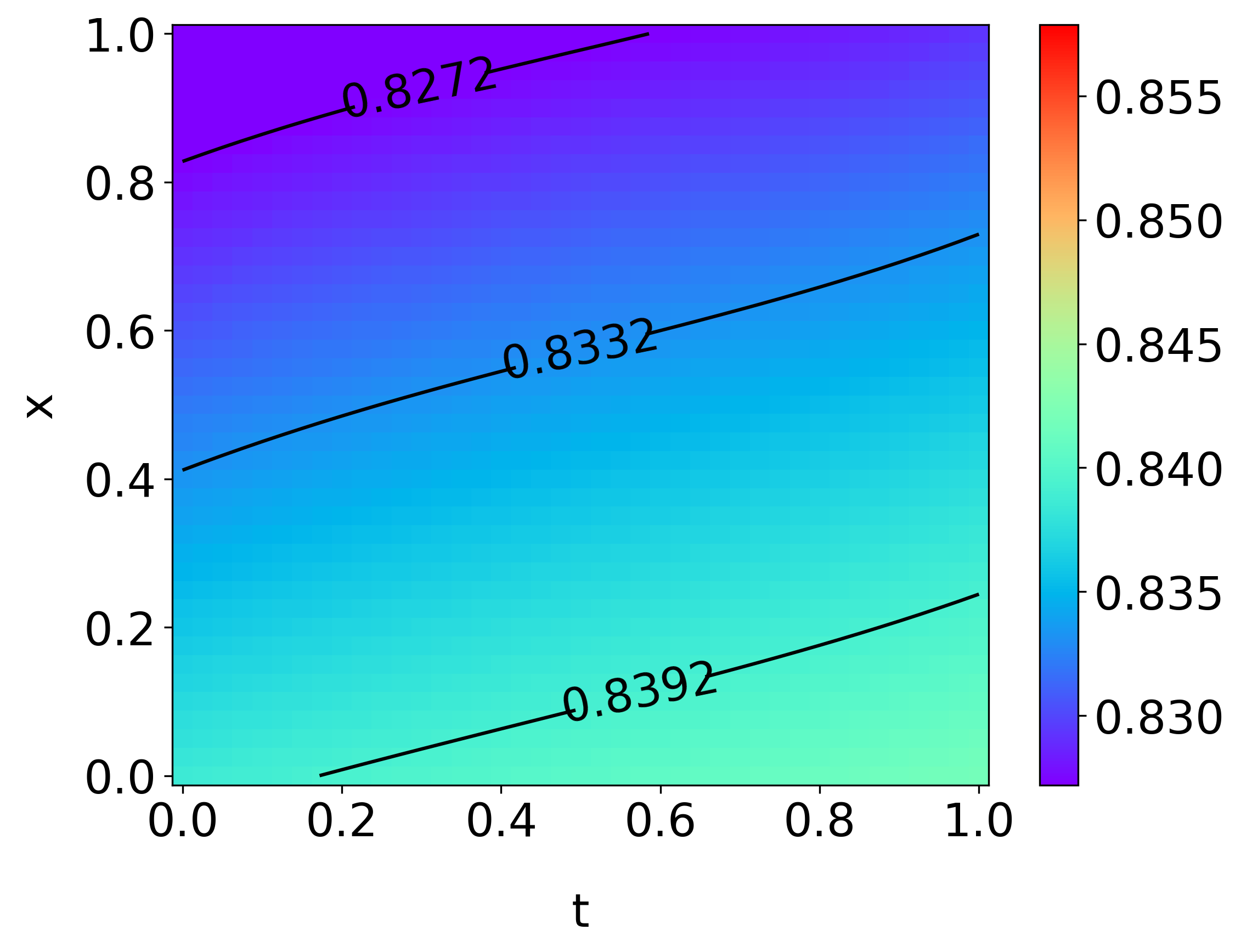}}

\caption{ Principal's value function (\emph{left}) relative error $V^P(t,x;0)/V^P(t,x)-1$ from \eqref{eq: HMsol}, and control (\emph{middle and right panels}) for the case study of Section \ref{subsec:HM}.
}
\label{fig:c=0 and c=0.5.}
\end{figure}

 For neural network surrogates of the value function $V(\cdot, \cdot)$ and the control $a(\cdot, \cdot)$, we use our DeepPAAC architecture described in Section~\ref{sec:NNA}. In the illustration below we take $L=3$ layers, with 32 neurons per layer. The networks $\theta_0^V, \theta_0^u$ are initialized using Glorot uniform initialization for the weights and zeros initialization for the biases.  The stopping rule is $\mathrm{Tol_{int}}=\mathrm{Tol_{ctrl}}=10^{-3}$.
For training we use batches of size $M=2000$, polynomial learning rate decay from $10^{-3}$ to $10^{-4}$ with power $0.8$. In the given example we need 6750 steps to satisfy the stopping criterion.

 In Figure \ref{fig:control for different c}, we plot the numerically obtained feedback control ${a}(t,x; \lambda)$ at a few locations  $(t,x)=({1\over 4},{1\over 4}),({1\over2},{1\over 2}),({3\over 4},{3\over 4})$ as we vary the mixture weight $\lambda$. The parameters in Figure \ref{fig:control for different c} are $\gamma_A=0.5$, $\gamma_P=1$ and $\tilde{\gamma}_P=2$. We observe that ${a}(t,x;\lambda) > a^{lin}(t,x;\lambda)$ while ${a}(t,x;0)$ and ${a}(t,x;1)$ match the analytical exact solution. The left panel of Figure \ref{fig:c=0 and c=0.5.} shows the relative error of $V^P$. The middle panel of Figure \ref{fig:c=0 and c=0.5.} shows ${a}(t,x; 0)$ as a function of $t,x$. As expected, the NN surrogate correctly learns that ${a}(t,x; 0)$ is (almost) a constant. The right panel of Figure \ref{fig:c=0 and c=0.5.} repeats this for ${a}(t,x; 0.5)$. We observe that when $\lambda=0.5$, ${a}(t,x; \lambda)$ increases in $t$ and decreases in $x$.

 \section{Multidimensional Case Studies}
\label{sec: Multidim-case}
\subsection{A Solely Continuous Payment}
\label{subsec: multi control-one dim spatial} 
We next consider a PA problem with a solely continuous payment which has been fully solved in \cite{Sung2025}. This example features both $\alpha_t,\beta_t$ as the Principal's controls, making it multi-dimensional unlike the toy example in Section \ref{subsec:HM}. Comparing to the generic case in Section \ref{sec:model}, we take  $b(t,a,x)=a, \sigma(t,x)=1$. We further assume no terminal lump sum payment $\xi \equiv 0$. The Agent's running utility is  
\begin{align}
U^A(t,\alpha_t,\beta_t,a_t,c_t)=k^A(c_t-{1\over 2}a_t^2),
\end{align}
where $k^A(\cdot)$ can be any increasing and concave function. The Principal's utility is linear
\begin{align}
U^P(t,\alpha_t,\beta_t,D_t)= D_t,
\end{align}
and terminal conditions are $\Phi^A(T,\xi,m^A)=m^A$ and $\Phi^P(T,\xi,X,m^P)=m^P$. From \cite{Sung2025}, we first solve the Agent's problem and get the optimal effort $\{a_t^*\}_{t\geq 0}$ and optimal consumption $\{c_t^*\}_{t\geq 0}$ to be: 
\begin{align*} \left\{ \begin{aligned}
&a_t^*=\beta_t+Z_t\\ 
&c_t^*={1\over2}(a_t^*)^2+K^A(1) 
\end{aligned}\right.
\end{align*}
where $K^A(x):=((k^A)')^{-1}(x)$. The process $(Z_t)$ is involved when we solve the Agent's problem. The following proposition also solves for the optimal control, which we will extend to a constrained case in Section 
\ref{sec:PA with constrained control}.
\begin{prop}
\label{prop: multidim-HJB}
The Principal's value function $V^P(t,w)$ satisfies the following HJB equation:
\begin{eqnarray}
\begin{cases}
\label{eq:Multidim-HJB}
V^P_t+\sup\limits_{\alpha,\beta,Z} \left\{(-C_0+{1\over2}Z^2-{1\over2}\beta^2-\alpha)V^P_w+{1\over2}Z^2V^P_{ww}+(1-\beta)(\beta+Z)-\alpha \right\} = 0,\\
V^P(T,w)=-w
\end{cases}
\end{eqnarray}
where $C_0=k^A(K^A(1))-K^A(1)$. The explicit solution is  $$V^P(t,w)=(C_0+{1\over 2})(T-t)-w.$$ The optimal Agent's control satisfies $a^*(t,w) = \beta^*(t,w)+Z^*(t,w)=1$.
\end{prop} 

To solve \eqref{eq:Multidim-HJB} numerically, we proceed to set up two different neural networks: one for approximating the value function $V^P,$ and the other for jointly modeling the control processes $\alpha,\beta,Z$ (an alternative, discussed below, is to train 3 separate networks for each of $\alpha, \beta, Z$).

For this example,  we have the vector convergence loss criterion
\begin{align} \label{eq:foc-3.3}
L_{ctrl}=\Big[-V^P_w - 1, \; -\beta V^P_w + 1 - 2\beta - Z, \; ZV^P_w + Z V^P_{ww} + 1 - \beta\Big].
\end{align}
For NN hyperparameters, we take training batches of $M=2000$, with $B= 10$ steps per epoch. The polynomial learning rate decays
from $10^{-3}$ to $10^{-4}$ with power $0.8$. The stopping rule is $\mathrm{Tol_{int}}=10^{-2}$ and $\mathrm{Tol_{ctrl}}=10^{-3}$, since the quantity $\| \bar{L}_{ctrl}( \theta^V_n, \theta^u_n; {\bar{t}}^\cdot, {\bar{w}}^\cdot) \|_{\infty}$ typically decays faster than the loss $\| \bar{L}_{int}( \theta^V_n; {\bar{t}}^\cdot, \bar{w}^\cdot) \|_{\infty}$ in our experiments.

\begin{figure}[htbp]
\centering
\subfigure[]{\includegraphics[width=6cm]{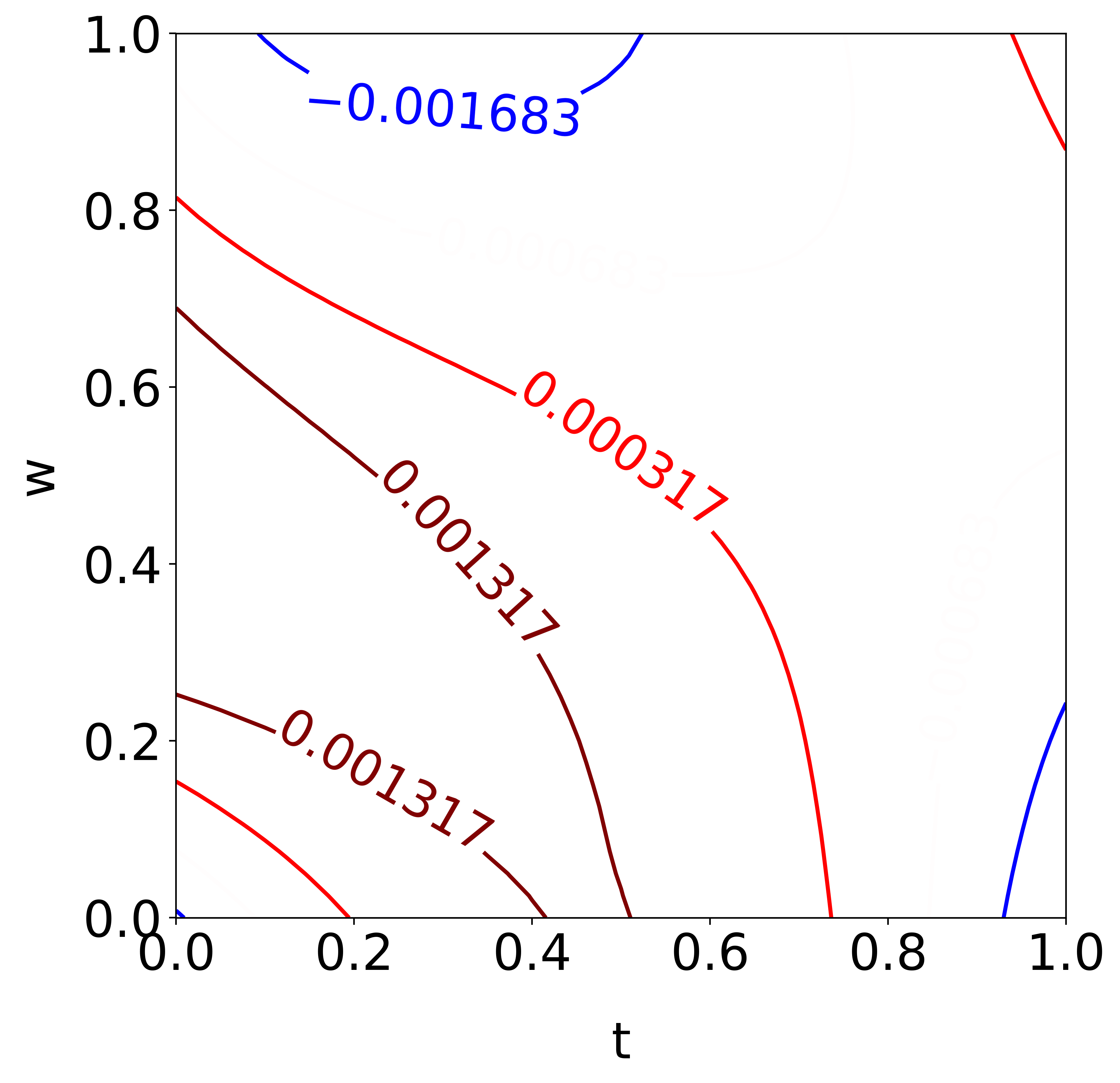}} 
\subfigure[]{\includegraphics[width=6cm]{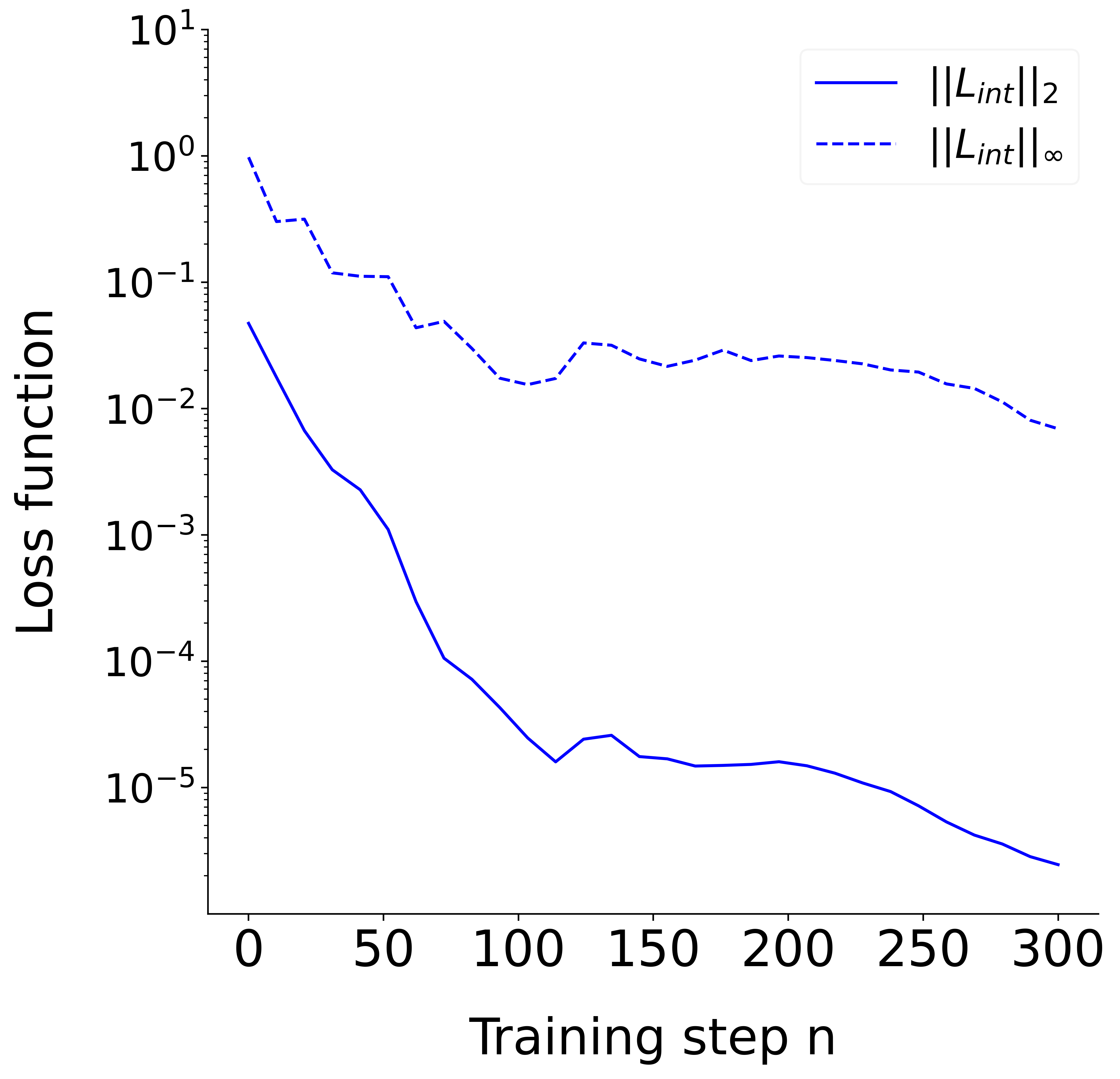}}
\\ 
\caption{\emph{Left (a):} Contour plot of the value function error metrics $\bar{L}_{int}(\theta_n^V;t,w)$ at $n=300$ training steps. \emph{Right (b):} $L_2$-norm and $L_\infty$-norm for $\bar{L}_{int}(\theta_n^V;\bar{t}^\cdot,\bar{w}^{\cdot})$ across $n$.} \label{fig: LH1 and LH2} 
\end{figure}

Choosing $C_0=0$, we solve \eqref{eq:Multidim-HJB}; numerical convergence to above tolerance is achieved in 300 training steps. Figure \ref{fig: LH1 and LH2} 
shows the pointwise  HJB residuals ${\bar{L}}_{int}(\theta_n^V;t,w)$ as a function of $(t,w)$ at $n=300$, as well as the $L_2$ and $L_\infty$-norm of the validation loss metric ${\bar{L}}_{int}(\theta_n^V;\bar{t}^{\cdot},\bar{w}^{\cdot})$  as a function of training steps $n$.

\subsection{DeepPAAC Variants}
\label{sec:DeepPAAC-Variants}

The defined DeepPAAC algorithm features many hyperparameters that can be finetuned during implementation. In this section, we describe several experiments that demonstrate the role of these hyperparameters.

First, we consider using DGM NN architecture and the relative performance of using multiple NN for learning $\bu$. To benchmark performance we compare the number of iterations needed to achieve convergence, based on the stopping rule of $\mathrm{Tol_{int}}=10^{-2}, \mathrm{Tol_{ctrl}}=10^{-3}$. We run each scheme five times (up to 10000 training steps each) in order to account for intrinsic variability of the SGD. For the DGM comparator,   we analogously take training batches of $M=2000$ in the interior, $M_T=1000$ for the terminal condition, $B=10$ SGD steps per epoch, and polynomial learning rate schedule for $\alpha_n$ from $10^{-3}$ to $10^{-4}$ with power $0.8$. 

\begin{figure}[htbp]
 \centering
\subfigure[$\| \bar{L}_{int}( \theta^V_n; {\bar{t}}^\cdot, \bar{w}^\cdot) \|_{\infty}$ across 5 runs]{\includegraphics[width=7cm]{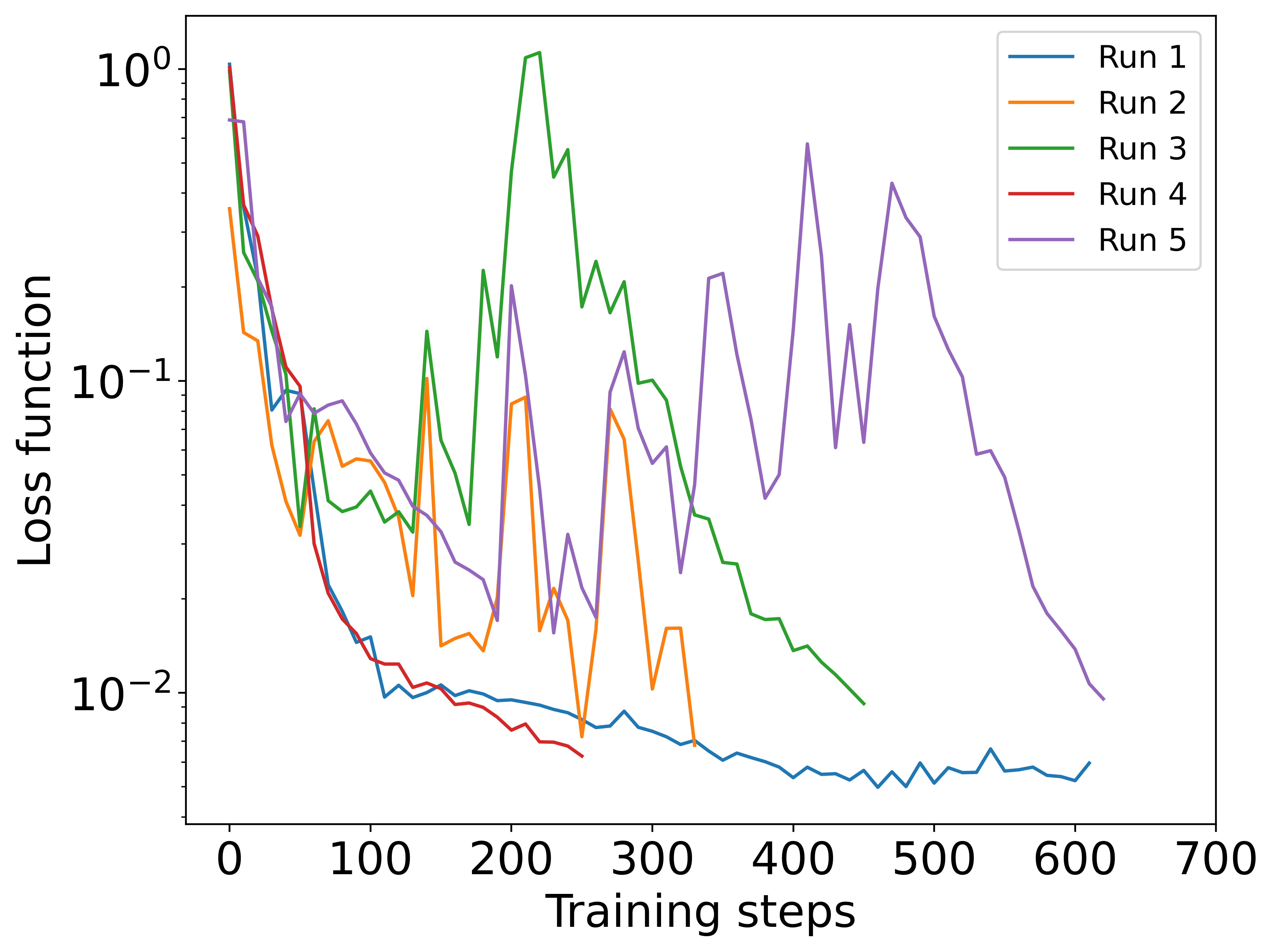}} 
\subfigure[$\| \bar{L}_{ctrl}( \theta^V_n, \theta^u_n; {\bar{t}}^\cdot, {\bar{w}}^\cdot) \|_{\infty}$ across 5 runs]{\includegraphics[width=7cm]{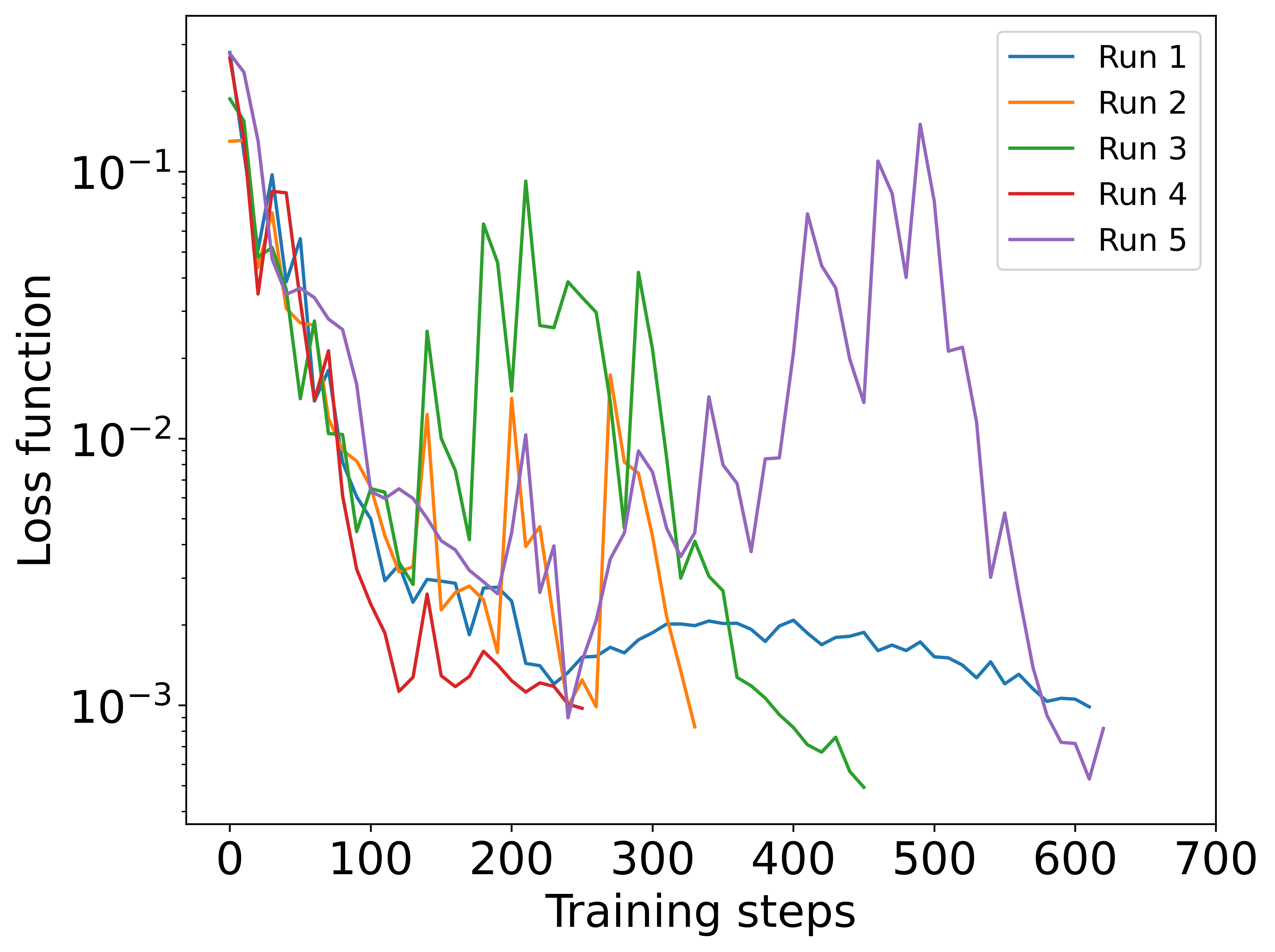}}
\caption{Validating convergence of the DeepPAAC loss functions. We show the $L_\infty$-norm for the value function HJB residual $\bar{L}_{int}( \theta^V_n; {\bar{t}}^\cdot, \bar{w}^\cdot)$ (\emph{left panel (a)}) and the control loss criterion $\bar{L}_{ctrl}( \theta^V_n, \theta^u_n; {\bar{t}}^\cdot, \bar{w}^\cdot)$ (\emph{right panel (b)}) as a function of step $n$, for the case study from Section \ref{subsec: multi control-one dim spatial}. Five runs of the DeepPAAC scheme with a single shared control NN.}\label{fig:test for NDGM}
\end{figure}

\begin{table}[ht]
\centering
\caption{Relative performance of different algorithms for case study of Section \ref{subsec: multi control-one dim spatial}. We report results across 5 runs to validate algorithm stability. The second column shows median number of SGD steps to achieve convergence and the third column the respective range across the 5 runs. Running time is in seconds. In the first row, one run never achieved the stopping threshold over 10,000 steps.  \label{tab:tuning}}
\begin{tabular}[t]{lrcr}
\toprule
Algorithm & Median \# of steps & Range of steps & Runtime \\
   & until convergence & across 5 runs & (seconds) \\
\midrule
DGM with shared control NN & $ 5770$ & $[1110,*] $ & 333  \\ 
DGM with three control NNs  & $ 4600$ & $[3740,7650]  $ & 248  \\
DeepPAAC  with shared control NN & $ 450$ & $[250,620] $ &16 \\
DeepPAAC  with three control NNs & $ 530$ & $[150,3270]$ & 22\\ \bottomrule
\end{tabular}
\end{table}
Table \ref{tab:tuning} shows that learning the three controls $\bu=(\alpha,\beta,Z)$ using a common shared NN performs better than utilizing three separate NNs in our DeepPAAC algorithm. Moreover, our algorithm converges faster than the DGM architecture. Figure \ref{fig:test for NDGM} plots the $L_\infty$-norm for $\bar{L}_{int}( \theta^V_n; {\bar{t}}^\cdot, \bar{w}^\cdot)$, and $\bar{L}_{ctrl}( \theta^V_n, \theta^u_n; {\bar{t}}^\cdot, \bar{w}^\cdot)$  of using DeepPAAC with a single shared NN for learning $\bu$. 

Next, we study the role of the sampling design. 

\begin{itemize}
    \item \textbf{Residual-based adaptive sampling (RAD)}: one strategy to speed up convergence is to preferentially train on samples where the loss functions are largest \cite{wu2023comprehensive}. To test this, at each epoch we generate $4 M$ candidate points in $(0,T) \times \cX$. After sorting them by the magnitude of the  PDE residual magnitude $L_1(t,w)$ we keep the top 20\% of the points (the ``high-error region'') which yield $0.8 M$ samples, as well as retain a random sub-sample of $0.2 M$ from the remaining ``low-error'' candidates. 
  
    \item \textbf{Non-uniform $t$-stratified sampling}: PDE residuals are generally larger for $t \ll T$, and hence it is more important to train on inputs where $t$ is small. To this end, we heuristically divide the time domain $t\in [0,1] = \cT_1 \cup \cT_2 \equiv  [0, 0.3] \cup  [0.3, 1]$. Half the samples are taken from $\cT_1$ and half from $\cT_2$, while the spatial domain remains uniformly sampled.
    
\end{itemize}

\begin{figure}[htbp]
\centering
\includegraphics[width=8cm,trim=0.2cm 0.7cm 0.2cm 0cm]{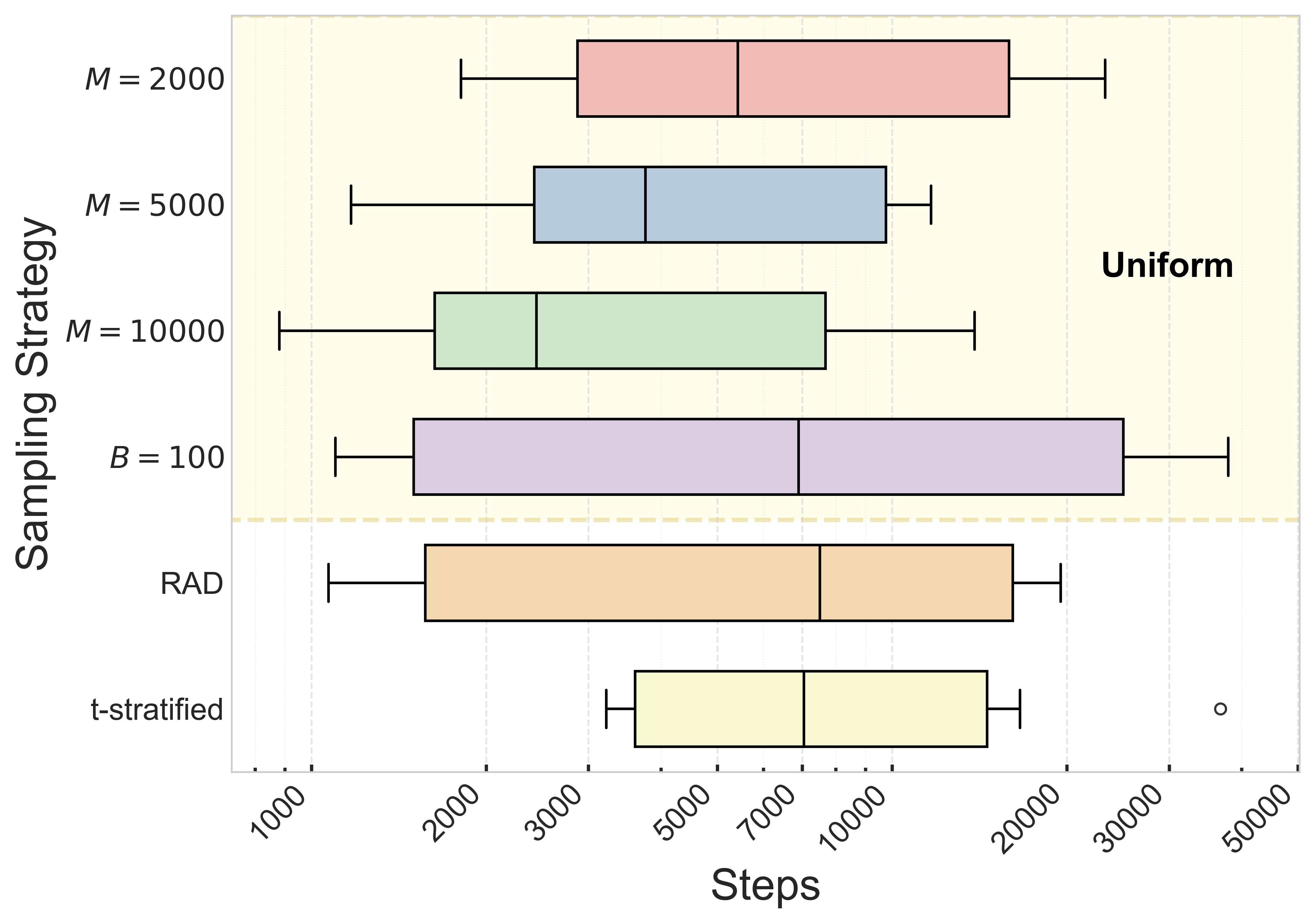}
\caption{Number of training steps needed to converge to loss thresholds $\mathrm{Tol_{int}}=\mathrm{Tol_{ctrl}}=10^{-3}$ under different training strategies for the case study of Section \ref{subsec: multi control-one dim spatial}. The $x$-axis is on the log-scale and the boxplots summarize performance across 10 runs of each variant. The top four variants use i.i.d.~Uniform sampling on $[0,T] \times \cX$, the bottom two (RAD and $t$-stratified) are described in text. \label{fig:sampling_comparison}}
\end{figure}

We also investigated the effect of the training set $M$, and the number $B$ of SGD steps taken at each epoch.  Under otherwise identical conditions, replacing only the sampling strategy led to the results shown in Figure \ref{fig:sampling_comparison}, with each strategy run 10 times. The stopping rule is $\mathrm{Tol_{int}}=\mathrm{Tol_{ctrl}}=10^{-3}$,  training batches of $M=2000$, ten steps per epoch $B=10$, polynomial learning rate decay from $10^{-3}$ to $10^{-4}$ with power $0.8$. For all experiments, the validation set used for the stopping criterion is fixed to $M_{V} = 2000$ with uniformly sampled points.

The parameter $M$ effectively acts as our batch size; we find that increasing the number of sampling points generally accelerates convergence, with only a marginal increase in the running time per training step. This happens because modern libraries distribute the loss function evaluation across $M$  training locations across the GPU/CPU memory. As the number of samples increases, the workload raises hardware utilization but there is minimal wall-clock impact until the memory resources are exhausted. At the same time, a large $M$ leads to a more accurate (i.e., less noisy) computation of the gradients and hence easier navigation of the loss landscape towards its minimum. As such, we recommend taking $M$ as big as possible given the computational capacity and memory resources of the hardware employed during training.

Finally, we also tried varying the number of gradient descent steps in each epoch. Experiments revealed that the variance across runs increased for larger $B$, such as $B=100$ rather than 10. Consequently, we recommend frequently resampling the training points.

\subsection{PA problem with Constrained Contracts}\label{sec:PA with constrained control}

Existing literature generally only considers unconstrained PA problems, since any constraints on Agent's remuneration immediately preclude any analytical solution. Nevertheless, realistic settings always feature some restrictions, both on what the Agent is able to do, and what the Principal is able to offer.  Thus we argue that constrained contracts are more realistic, but have been viewed as intractable due to lack of numerical methods that our DeepPAAC scheme overcomes, opening the door to many new versions of contracting. The implementation for this section is publicly available at \url{https://github.com/Reedcgx/DeepPAAC}.

In this section, we consider a problem similar to Section \ref{subsec: multi control-one dim spatial}. However, we add constraints on $\beta_t$ and $\beta_t+Z_t$ by restricting them to satisfy
\begin{equation}
\underline{\beta}\leq \beta_t \leq \overline{\beta}, \qquad \underline{u} \leq \beta_t+Z_t\leq \overline{u}. \label{eq:constraint}
\end{equation}
Recall from Proposition \ref{prop: multidim-HJB} that the Agent's optimal effort is $a^*_t= \beta_t + Z_t$, so constraining the sum $\beta(t,w)+Z(t,w)$ is equivalent to the Agent having an a priori bounded effort in the range $[\underline{u}, \bar{u}]$. Similarly, constraining $\beta$ corresponds to limiting how much payment is given based on the observable output $(X_t)$, cf.~\eqref{eq:St-dynamics}. Since there are three variables in the Principal's problem $\alpha,\beta,Z$,  the above constraints are equivalent to keeping a single 
unbounded state and bounding the other two. 

Similar to Proposition \ref{prop: multidim-HJB} in Section \ref{subsec: multi control-one dim spatial}, we have the following explicit solution:

\begin{prop}
Under constraints \eqref{eq:constraint},
the Principal's value function $V^P(t,w)$ satisfies the HJB equation:
\begin{eqnarray}  \hspace*{20pt}
\begin{cases}
\label{eq:Multidim-HJBwithconstraint}
V^P_t+\sup\limits_{\substack{\alpha,\ \underline{\beta}\leq \beta\leq \bar{\beta}, \\ \underline{u}\leq \beta+Z\leq \overline{u}}} \left\{(-C_0+{1\over2}Z^2-{1\over2}\beta^2-\alpha)V^P_w+{1\over2}Z^2V^P_{ww}+(1-\beta)(\beta+Z)-\alpha \right\} =0\cr
V^P(T,w)=-w,
\end{cases}
\end{eqnarray}
where $C_0=k^A(K^A(1))-K^A(1)$ and $K^A(x):=((k^A)')^{-1}(x)$.  The explicit solution of \eqref{eq:Multidim-HJBwithconstraint} is: 

(1) Nonbinding constraints: If $\overline{u}\geq 1$ and $\underline{u}\leq 1$, then 
\begin{align} \left\{ \begin{aligned}
& V^P(t,w) =(C_0+{1\over 2})(T-t)-w \\  
& \beta^*(t,w)+Z^*(t,w)\equiv 1. \end{aligned} \right.\end{align}

(2) Binding constraints: 

(i) If $\overline{u}< 1$, then 
\begin{align}\left\{ \begin{aligned}
&V^P(t,w) = \left[C_0-{1\over 2}(\overline{u})^2+\overline{u} \right](T-t)-w \\ 
&\beta^*(t,w)+Z^*(t,w)  \equiv \overline{u}. 
\end{aligned}\right. \end{align}

(ii) If $\underline{u}>1$, then 
\begin{align}\left\{ \begin{aligned}
&V^P(t,w) = \left[C_0-{1\over 2}(\underline{u})^2+\underline{u} \right](T-t)-w \\ 
&\beta^*(t,w)+Z^*(t,w)  \equiv \underline{u}. 
\end{aligned}\right. \end{align}
\end{prop} 

\proof The terminal condition $V^P(T,w)=-w$ can be derived similarly as in \cite{Sung2025}. Observe that in the HJB equation (\ref{eq:Multidim-HJBwithconstraint}), since there is no constraint on $\alpha$, we must have $V^P_w(t,w)\equiv -1$ and $V^P_{ww}(t,w) \equiv 0$. Therefore, we only need to compute
\begin{align*}
    &\sup\limits_{\underline{\beta}\leq \beta\leq \overline{\beta},\underline{u}\leq \beta+Z\leq \overline{u}} \left(C_0-{1\over2}Z^2+{1\over2}\beta^2+(1-\beta)(\beta+Z)\right) \\ 
    =&\sup\limits_{\underline{\beta}\leq \beta\leq \overline{\beta},\underline{u}\leq \beta+Z\leq \overline{u}}
    \left(C_0+{1\over 2}-{1\over 2}(\beta+Z-1)^2 \right).
\end{align*}
 The above operator is a quadratic function of $\beta+Z$, and the optimal value can be obtained correspondingly, depending on whether  the interval $[\underline{u},\overline{u}]$ contains $1$ or not. The rest of the steps to find the solution $V^P(t,w)$ are straightforward. \qed

In order to numerically handle constraints on PA controls, we extend the  PIA algorithm in Section \ref{subsec:algo}. Namely, we augment the loss function of the control NN to additionally penalize (pointwise at $(t,\bx)$) controls that violate the constraints. To this end we introduce a penalty function $P(\bu)$
that is zero when $\bu \in \mathcal{U}(t,\bx)$ and is strictly positive otherwise. Then we add the penalties due to not satisfying the constraint over the training batch to the loss metric in Step 5, which encourages taking SGD steps towards the feasible region where all the constraints hold.
The resulting loss metric for the control NN $\bu(\cdot,\cdot; \theta^u_n)$ is 
\begin{align}
\label{eq:constraint-penalty}
    L_u(\theta_n^u)=&-\sum_{m=1}^{M_V}  \Bigl[\cL^{\bu(t_n^m,\bx_n^m;\theta_n^u)}v(t_n^m,\bx_n^m;\theta_n^V)+ F(t_n^m,\bx_n^m,\bu(t_n^m,\bx_n^m;\theta_n^u))\Bigr] \nonumber  \\ 
    &+\lambda_{Pen} P(\bu(t_n^m,\bx_n^m;\theta_n^u)), 
\end{align}
for some penalty parameter $\lambda_{Pen}$, taken to be 1 by default.

In our example above, since the constraints are on $\beta$ and $Z$, we have the convergence loss criterion reducing to
\begin{align} \label{eq:foc-4.3}
 &  L_{ctrl}=\Big[-V^P_w - 1\Big].
\end{align}

\begin{figure}[!ht]
\centering
\subfigure[${\beta}(t,w)$]{\includegraphics[width=0.28\textwidth]{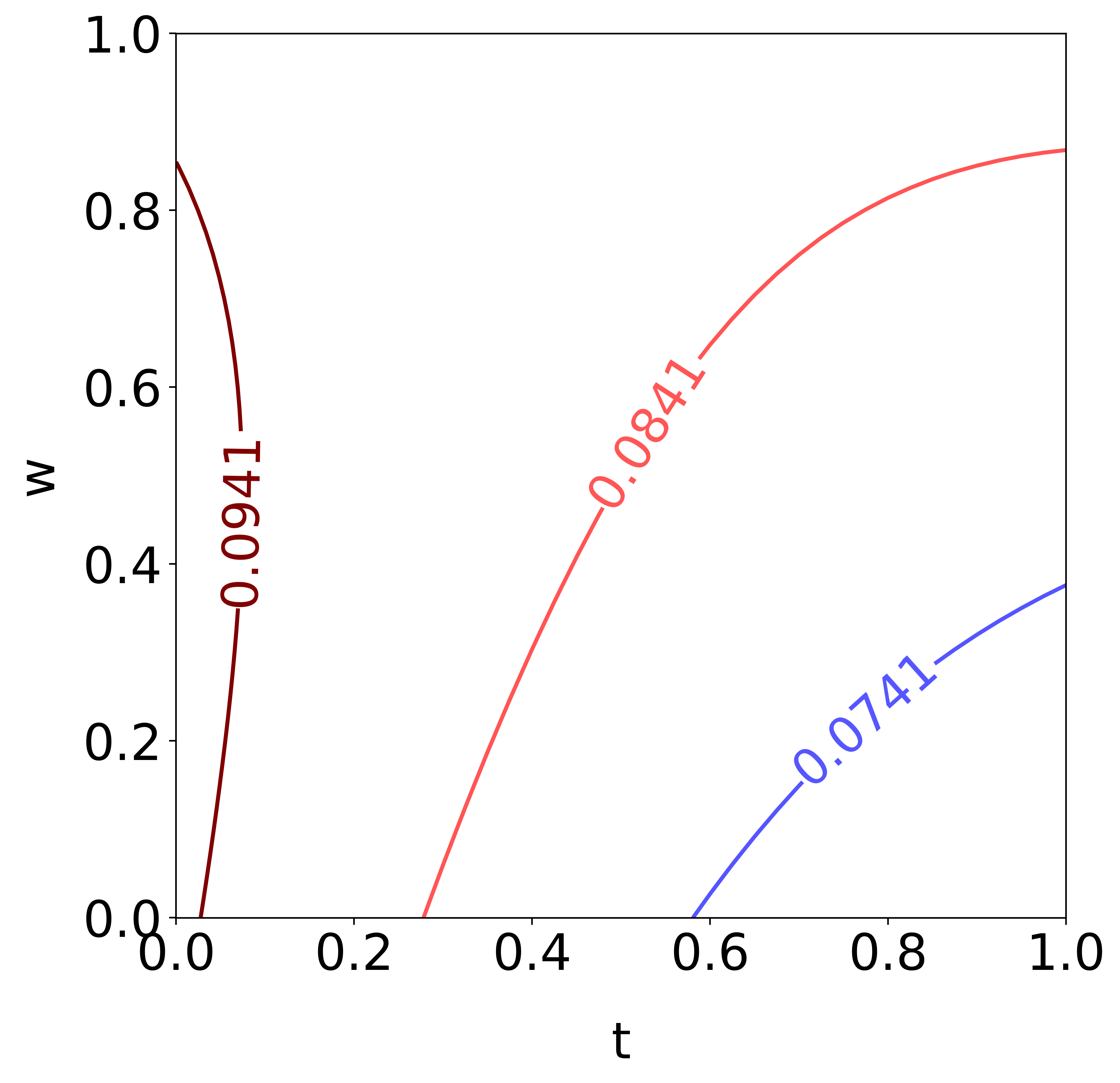}} 
\subfigure[${\beta}(t,w)+{Z}(t,w)$]{\includegraphics[width=0.28\textwidth]{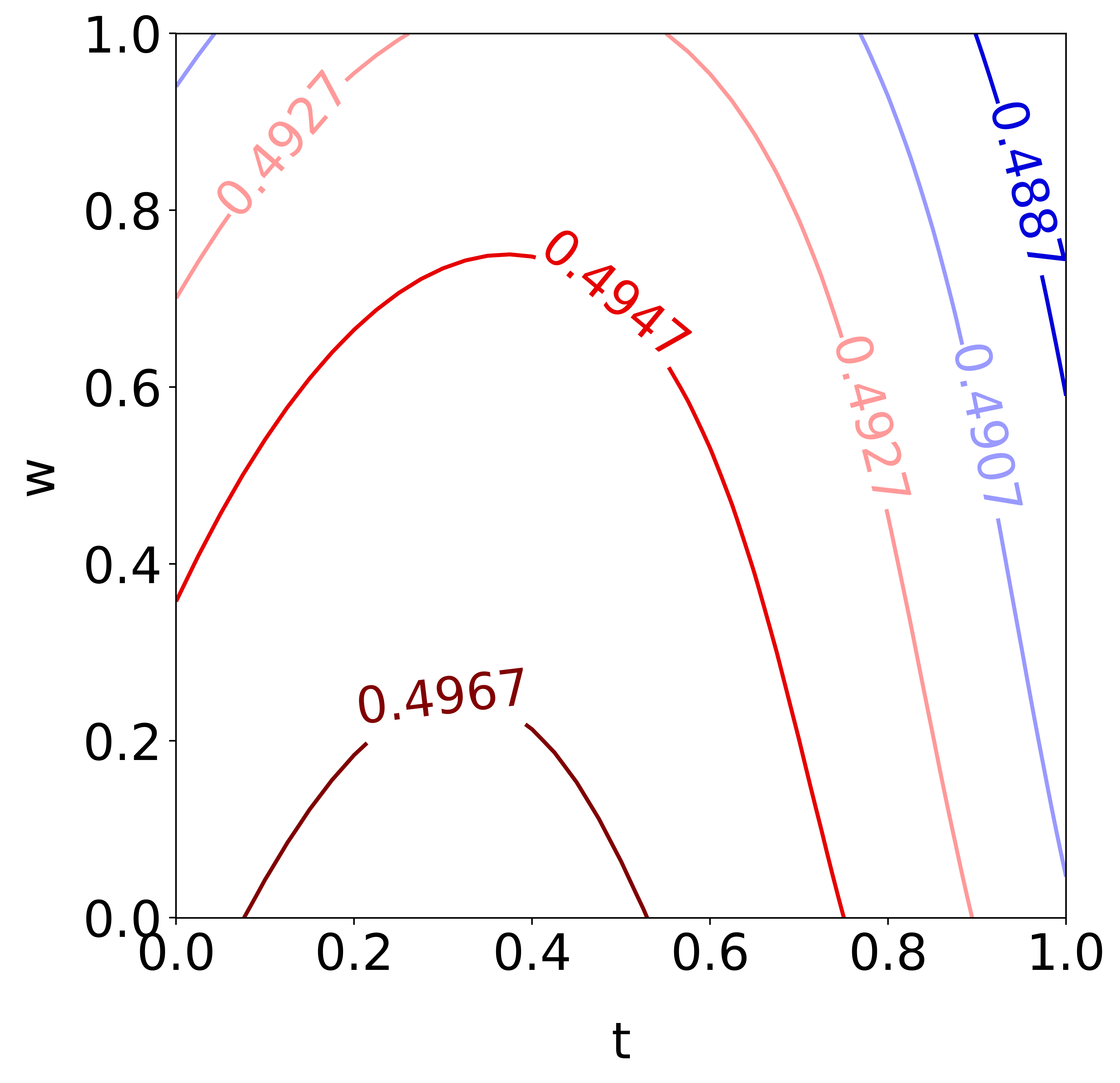}}
\subfigure[Loss function]{\includegraphics[width=0.28\textwidth]{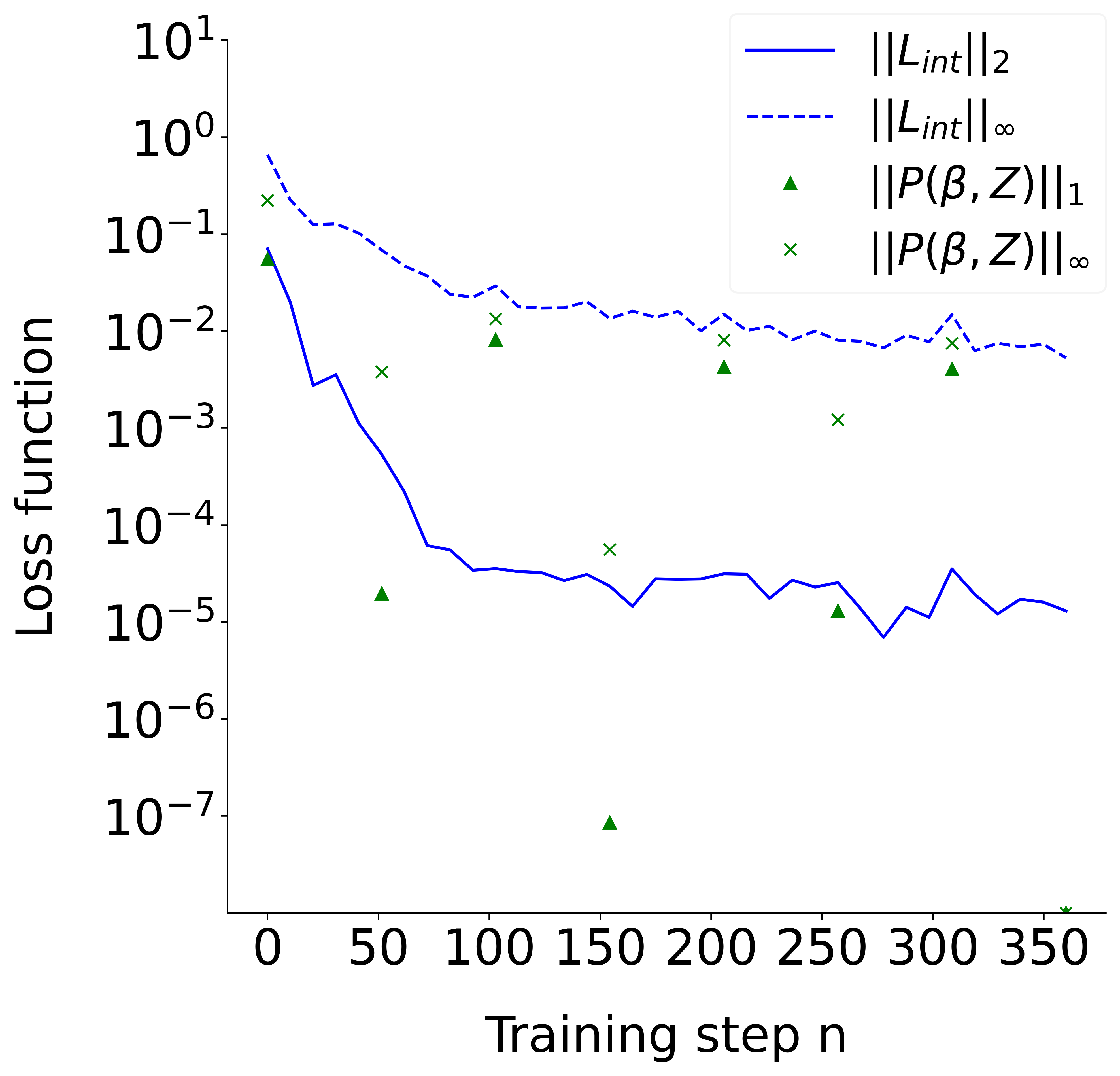}}\\
\subfigure[${\beta}(t,w)$]{\includegraphics[width=0.28\textwidth]{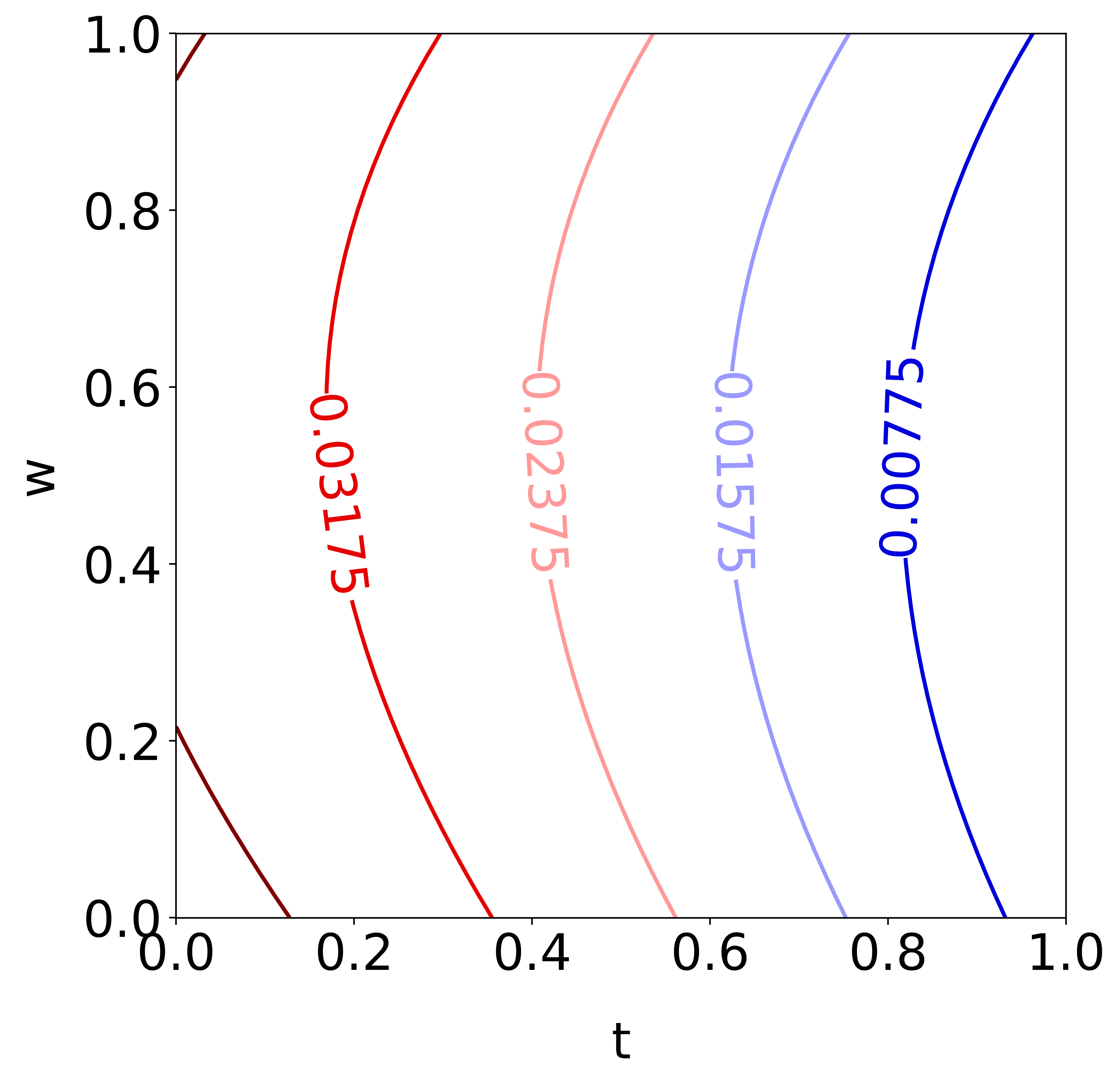}} 
\subfigure[${\beta}(t,w)+{Z}(t,w)$]{\includegraphics[width=0.28\textwidth]{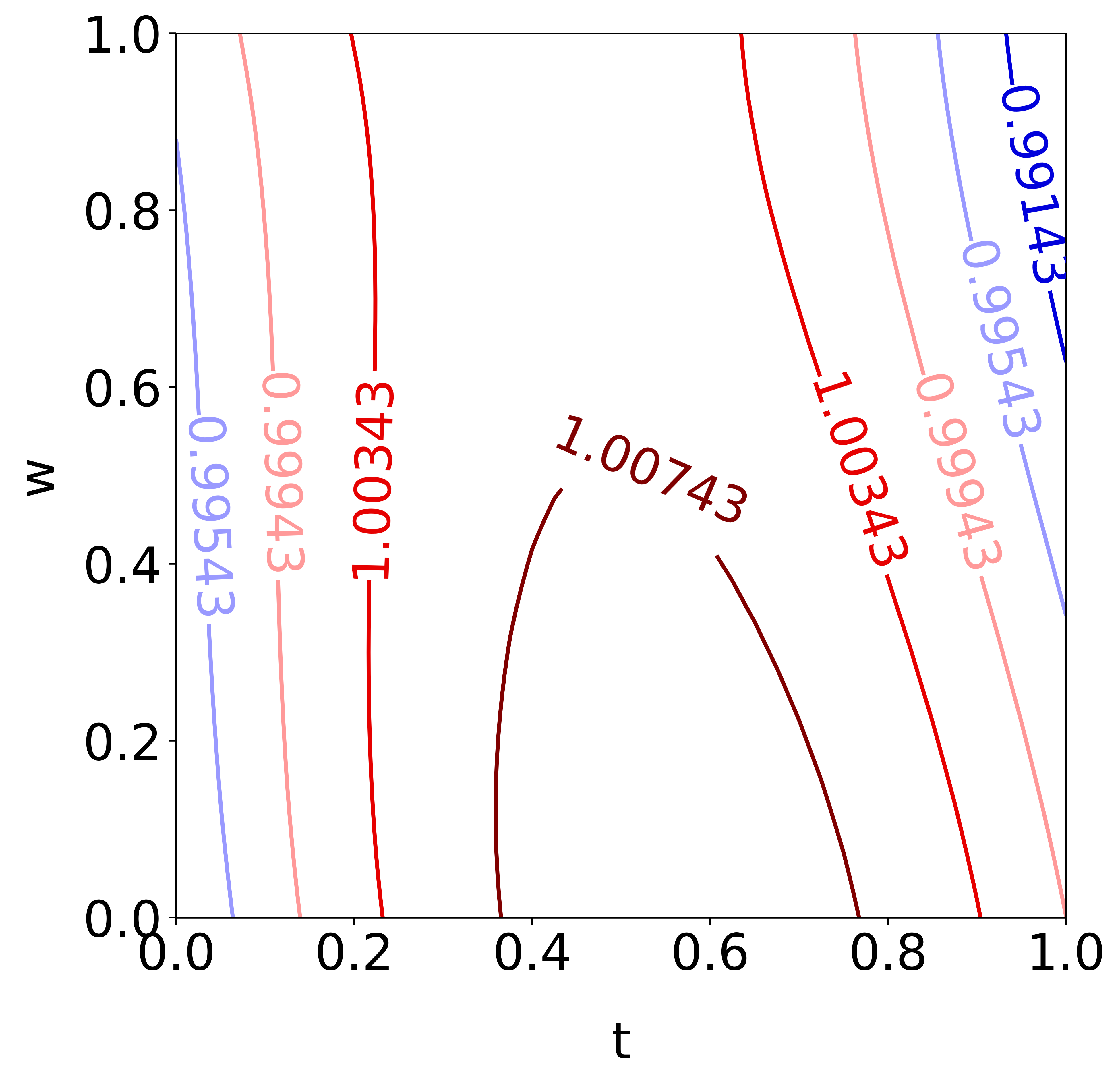}}
\subfigure[Loss function]{\includegraphics[width=0.28\textwidth]{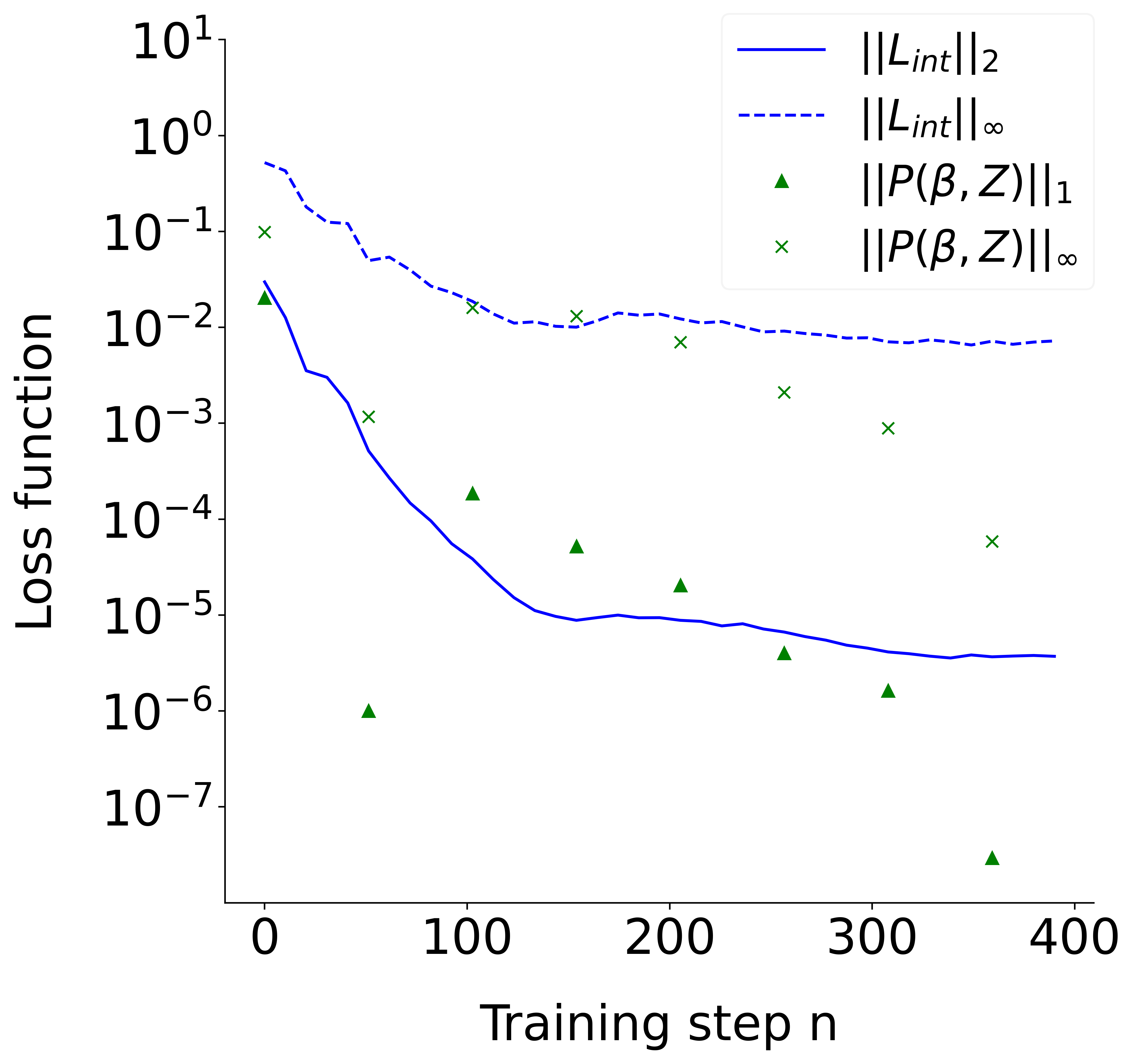}} \\
\subfigure[${\beta}(t,w)$]{\includegraphics[width=0.28\textwidth]{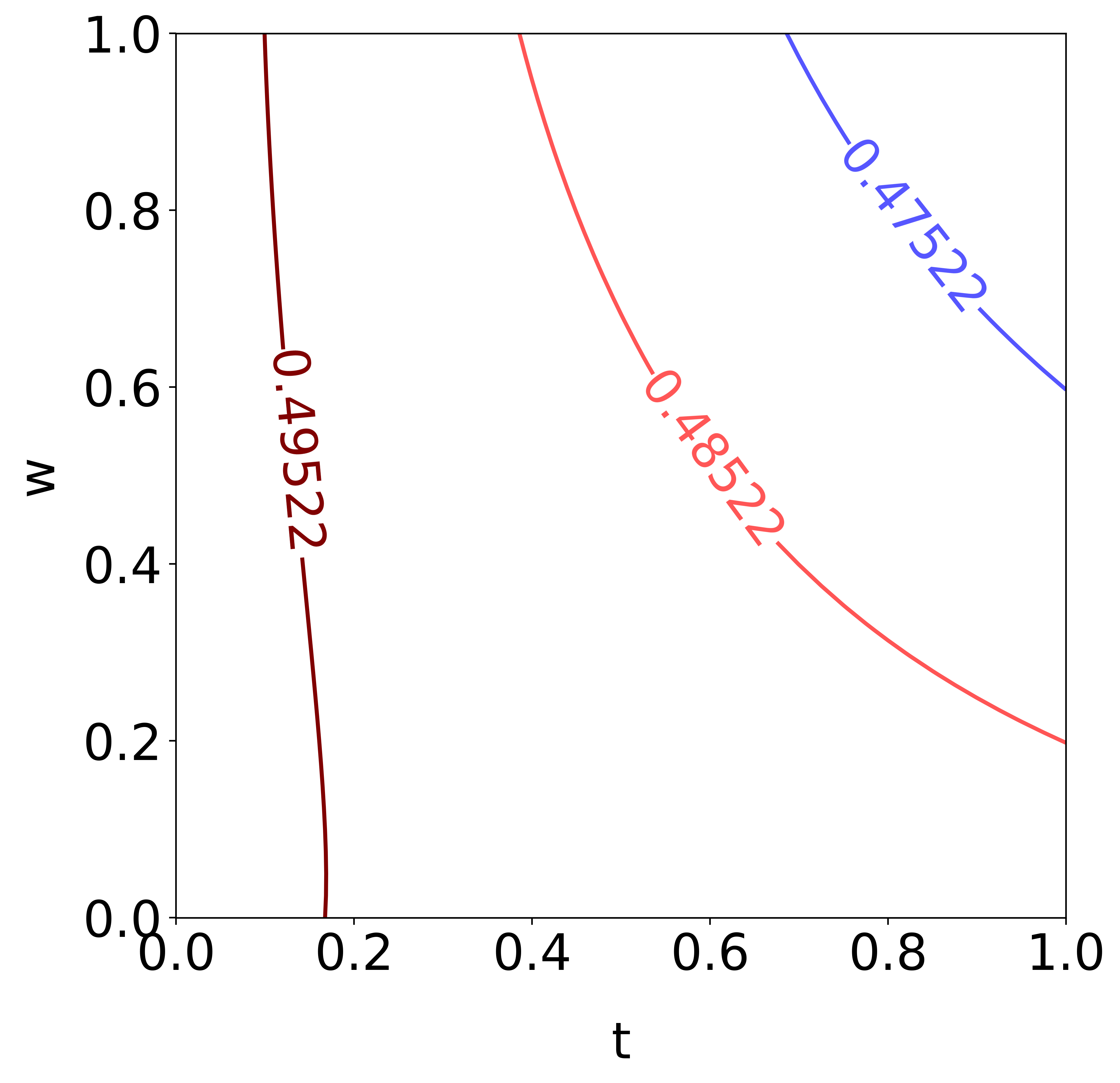}} 
\subfigure[${\beta}(t,w)+{Z}(t,w)$]{\includegraphics[width=0.28\textwidth]{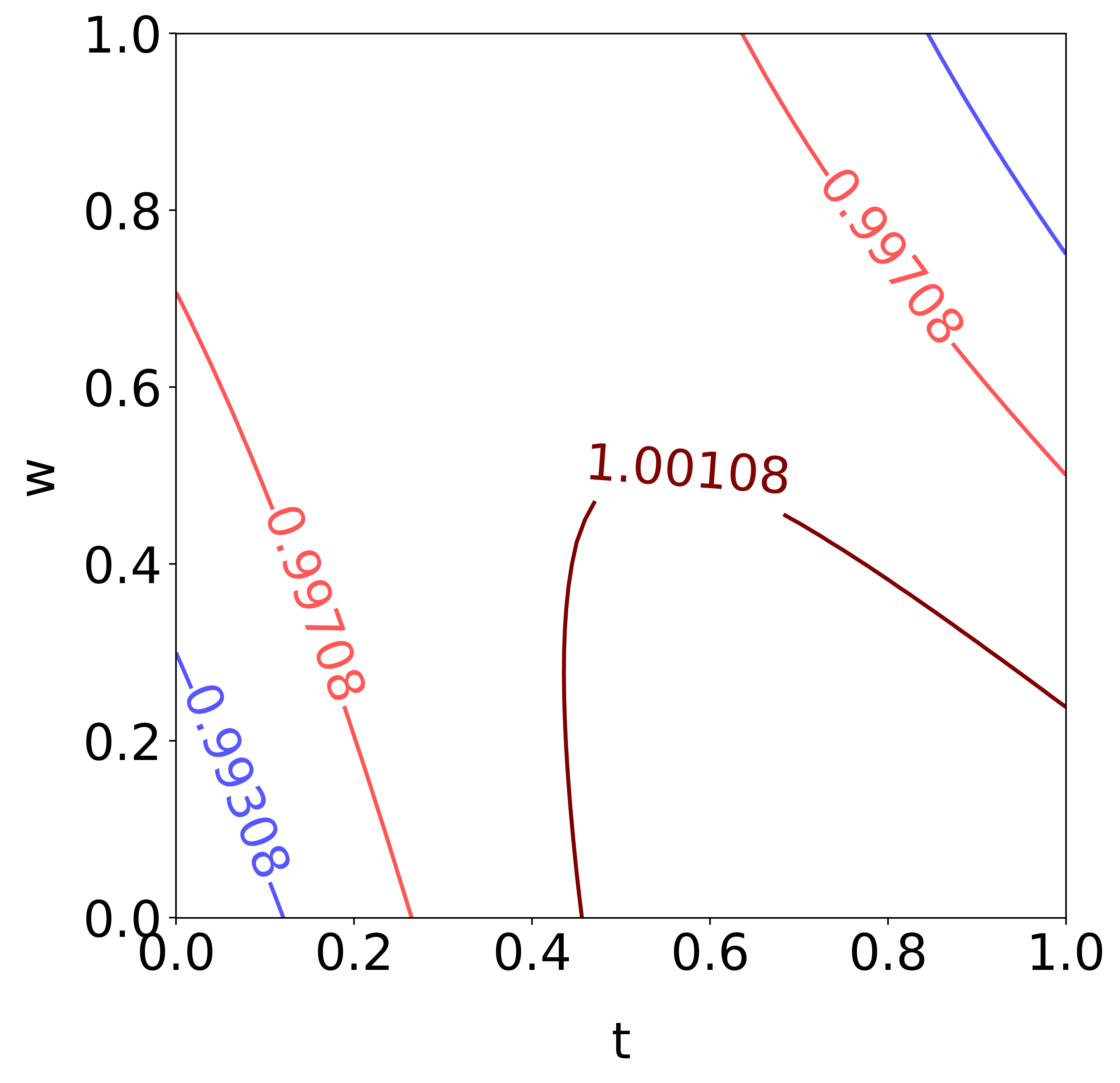}}
\subfigure[Loss function]{\includegraphics[width=0.28\textwidth]{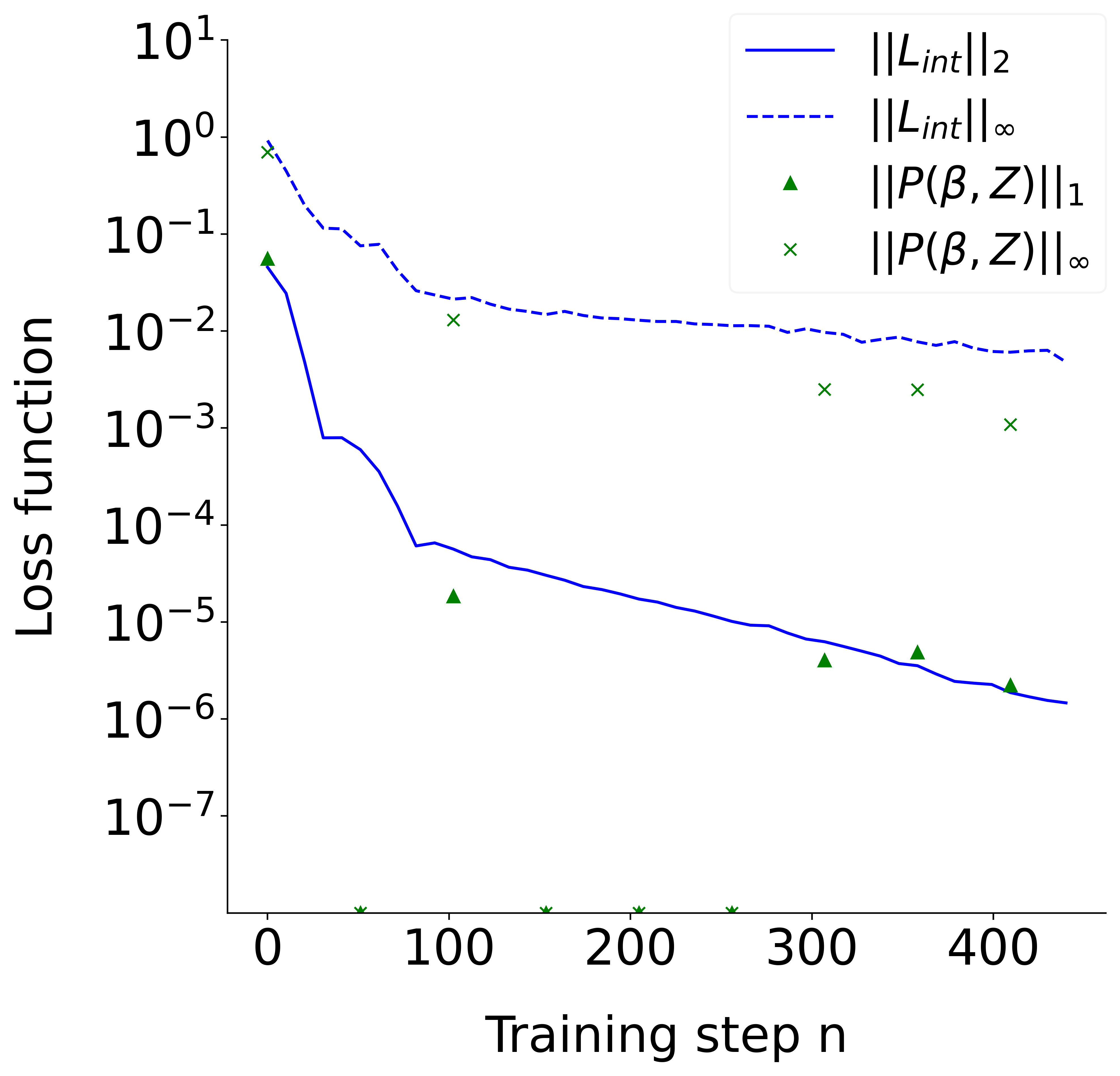}}
\caption{Validating the solution of the constrained PA case study from Section \ref{sec:PA with constrained control} across three different constraint settings. \emph{Top row:} $0\leq \beta\leq 0.1, 0\leq \beta+Z\leq 0.5$. \emph{Middle row:} $0\leq \beta\leq 0.2$. \emph{Bottom row:} $\beta \leq 0.5, \beta+Z\leq 1.2$. The right-most panels $(c)-(f)-(i)$  trace the $L_2$-norms and $L_\infty$-norms of $\bar{L}_{int}(\theta_n^V;\bar{t}^\cdot,\bar{w}^\cdot)$  as well as the penalty functions $\overline{P}_1$, $\overline{P}_{\infty}$. Numerical convergence is achieved in $350/380/430$ training steps for cases 1-3 respectively. \label{fig:case1}}
\end{figure}

To solve \eqref{eq:Multidim-HJBwithconstraint} numerically,  we use the architecture design as in Section \ref{subsec: multi control-one dim spatial}. The stopping rule is $\mathrm{Tol_{int}}=10^{-2}, \mathrm{Tol_{ctrl}}=10^{-3}$ and $\overline{P}_{\infty}=0$  where the latter is the norm of the penalty function on the validation set $\{\bar{t}^m,\bar{w}^m\}_{m=1}^{M_V}$:
\begin{align}\nonumber
\overline{P}_1 := \frac{1}{M_V} \sum_{m=1}^{M_{V}} {P}(\bar{t}^m,\bar{w}^m,u(\bar{t}^m,\bar{w}^m;\theta_n^u)), \qquad \overline{P}_\infty := \max_{m \le M_V} {P}(\bar{t}^m,\bar{w}^m,u(\bar{t}^m,\bar{w}^m;\theta_n^u)).
\end{align}
We use $M=2000 $ training points and $B=10$ steps per epoch. Learning rate schedule is with polynomial decay
from $10^{-3}$ to $10^{-4}$ with power $0.8$. 
In the following three numerical examples, we illustrate several different constraint combinations. Throughout  we choose $C_0=0$.

\textbf{Case 1:} Take $\underline{\beta}=0,\overline{\beta}=0.1,\underline{u}=0,\overline{u}=0.5$ with the penalty function:
\begin{align}\nonumber
P(\beta,Z)=(0-\beta)^+ +(\beta-0.1)^+ +(0-\beta-Z)^+ +(\beta+Z-0.5)^+.
\end{align}
In top row of Figure \ref{fig:case1},  we plot the numerically obtained optimal control ${\beta}(t,w)$, as well as the sum ${\beta}(t,w)+{Z}(t,w)$ (so one can read off ${Z}$ as well). We observe that the constraint on $\beta+Z$ is binding, ${\beta}(t,w)+{Z}(t,w)\simeq \overline{u}=0.5$ as expected. But the constraint of ${\beta}(t,w)$ is not binding. Panel (c) shows the $L_1$ and $L_\infty$ norms of the penalty function $\overline{P}_1,\overline{P}_\infty$
as a function of steps $n$ until convergence.

\textbf{Case 2}: Take $\underline{\beta}=0,\bar{\beta}=0.2$ and no constraint on Z. The penalty function is 
\begin{align*}
P(\beta,Z)=(0-\beta)^+ + (\beta-0.2)^+.
\end{align*}
The middle row of Figure \ref{fig:case1} shows that the constraint is no longer binding in the interior, and ${\beta}(t,w)$ can be as low as 0.16 for large $t$. 

\textbf{Case 3:} Take $\bar{\beta}=0.5,\overline{u}=1.2$ with the penalty function
\begin{align*}
P(\beta,Z)=(\beta-0.5)^+ + (\beta+Z-1.2)^+.
\end{align*}
The bottom row of Figure \ref{fig:case1} shows that (as expected) the constraints on $\beta$ and $\beta+Z$ do not bind, and the sum ${\beta}(t,w) + Z(t,w)$ remains 1 everywhere.

Penalty function scaling: in \eqref{eq:constraint-penalty} the penalty acts as a third additive term for the loss criterion, raising the question of selecting the scaling parameter $\lambda_{Pen}$. Table \ref{tab:penalty} below compares across three choices of $\lambda_{Pen} \in \{0.2, 1, 5\}$. We observe that different scalings work better for different cases; when the constraint is strong, it is beneficial to take a larger $\lambda_{Pen}$, otherwise it is better to concentrate on the other terms of $L_u(\cdot)$. In offline experiments not reported directly here, we also tried other penalty specifications, such as $P(\beta)=\exp( (\beta-\underline{\beta})^+) $ or $P(\beta)=((\beta-\underline{\beta})^+)^2$ that more strongly penalize large violations; however these modifications do not work empirically as well as the basic linear penalty. 

\begin{table}[ht]
\centering
\caption{Relative performance of different penalty functions for the constrained cases.
We report results across 5 runs to evaluate stability.
The third column shows the median number of SGD steps until convergence, and the fourth column shows the corresponding range across the 5 runs.
Each asterisk (*) in the range column denotes one run that reached the maximum step limit of 1500 without satisfying the stopping criterion. 
\label{tab:penalty}}
\begin{tabular}{llrr}
\toprule
Case & Penalty & Median \# of steps & Range of steps \\
 & & until convergence & across 5 runs \\
\midrule
\multirow{3}{*}{Case 1}
& \textbf{Base} $\lambda_{Pen}=1$ & 280 & [130, 960] \\
& $\lambda_{Pen}=0.2$ & 1500 & [1500, *****] \\
& $\lambda_{Pen}=5$ & 800 & [430, 1080] \\
\midrule
\multirow{3}{*}{Case 2}
& \textbf{Base} & 390 & [220, 460] \\
& 0.2$\times$ Penalty & 390 & [320, 1430] \\
& 5$\times$ Penalty & 390 & [170, 740] \\
\midrule
\multirow{3}{*}{Case 3}
& Base & 550 & [160, *] \\
& 0.2$\times$ Penalty & 520 & [380, 1010] \\
& \textbf{5$\times$ Penalty} & 210 & [110, 1020] \\
\bottomrule
\end{tabular}
\end{table}

\subsection{Two-dimensional Problem with Scalar Control and Explicit Solution} \label{sec:mul}

In this subsection, we test our algorithm on Example 5.5 in Cvitani{\'c}, Possama{\"\i} and Touzi \cite{cvitanic2018dynamic}, which is a two-dimensional Principal-Agent problem with explicit solution. The implementation for this section is publicly available at \url{https://github.com/Reedcgx/DeepPAAC}. Translating into our notation from Section \ref{sec:model}, we take  $b(t,a,x)=\sigma ax, \sigma(t,x)=\sigma x$ 
which corresponds to Geometric Brownian motion dynamics for $\{X_t\}$, quadratic running cost for the Agent, and zero running cost for Principal:
$$U^A(t,\alpha,\beta,a,c)=-{1\over 2}a^2, \qquad   U^P(t,\alpha,\beta,D)=0.$$
For the terminal lump-sum payment, there is a risk-averse Agent utility function \\ $\Phi^A(T,\xi,m^A)= \log \xi$ and risk-neutral Principal's utility function $\Phi^P(T,\xi,X,m^P)=X-\xi$.

The Agent's optimal wealth process $\{W_t\}$ satisfies
$$dW_t={1\over 2}\sigma^2 X_t^2Z_t^2 \,dt+\sigma X_t Z_t \ dB_t^a, \qquad W_T  = \log \xi.$$

\noindent It is shown in \cite{cvitanic2018dynamic} that the Agent's optimal control is
$a^*_t=\sigma X_t Z_t$ and the HJB equation for the Principal's value function $V^P(t,x,w)$ is
\begin{eqnarray}
\begin{cases}
\label{eq:CPT example-explicit}
V^P_t+\sigma^2 x^2\sup\limits_{Z \in \mathbb{R}} \left\{ZV^P_x+{1\over 2}Z^2V^P_w+{1\over 2}(V^P_{xx}+Z^2V^P_{ww})+ZV^P_{xw}\right\} =0,\\
V^P(T,x,w)=x-e^w. 
\end{cases}
\end{eqnarray}
It is reported in \cite{cvitanic2018dynamic} that \eqref{eq:CPT example-explicit} yields an explicit solution:
$$V^P(t,x,w)=x-e^w+{1\over 4}e^{-w}x^2\big(e^{\sigma^2(T-t)}-1\big), \qquad Z^*(t,x,w)={1\over 2}e^{-w}.$$
For this example, the convergence loss criterion is the gradient in $Z$ of the Hamiltonian:
$$ L_{ctrl}=\Big[(V^P_w+V^P_{ww})Z+V^P_x+V^P_{xw}\Big].$$

Below we work with $\sigma=1$.

\begin{figure}[!htb]
\centering
\subfigure[absolute error of DeepPAAC value function ]{\includegraphics[width=0.45\textwidth]{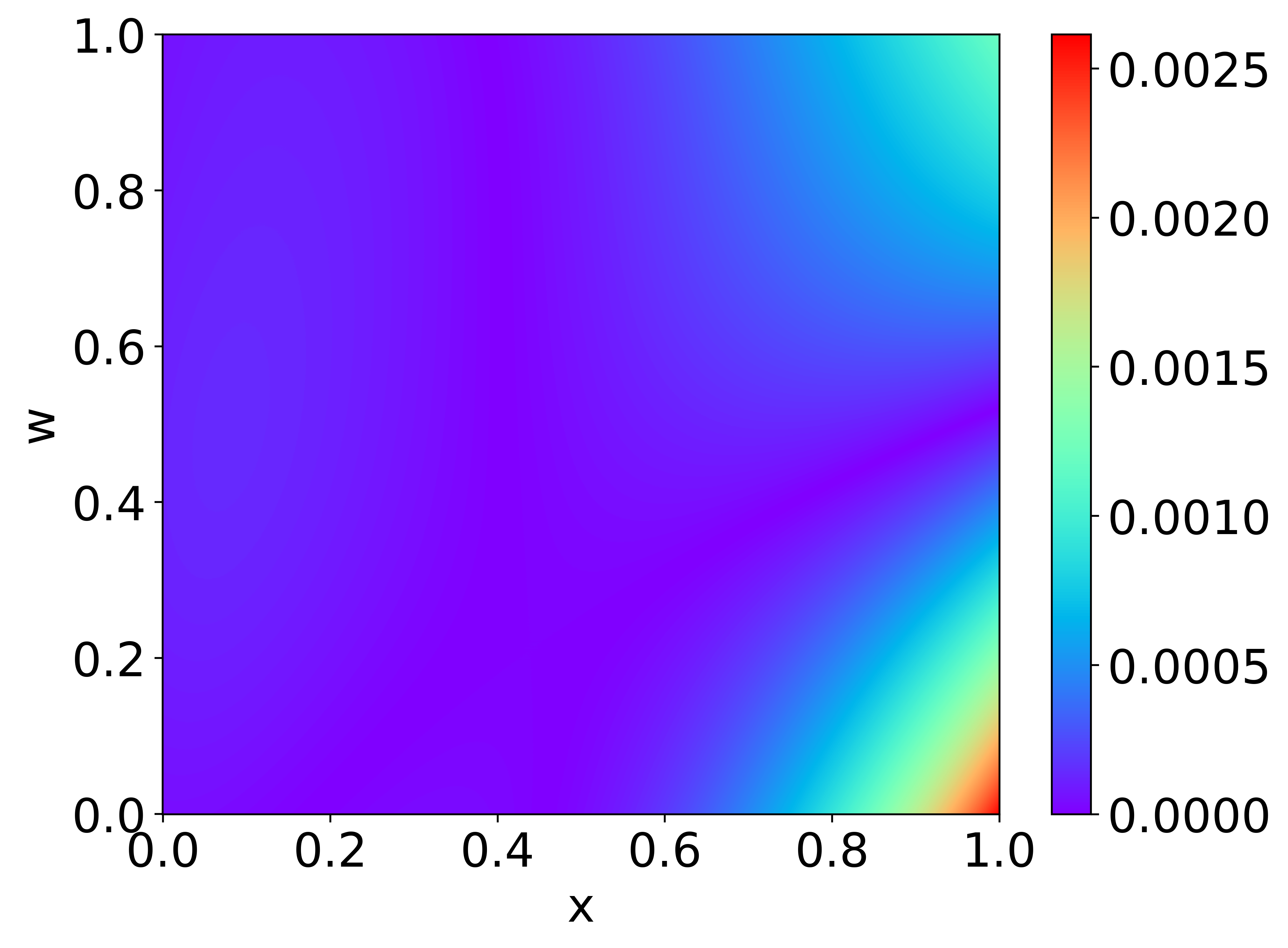}} 
\subfigure[absolute error of DeepPAAC optimal control ]{\includegraphics[width=0.45\textwidth]{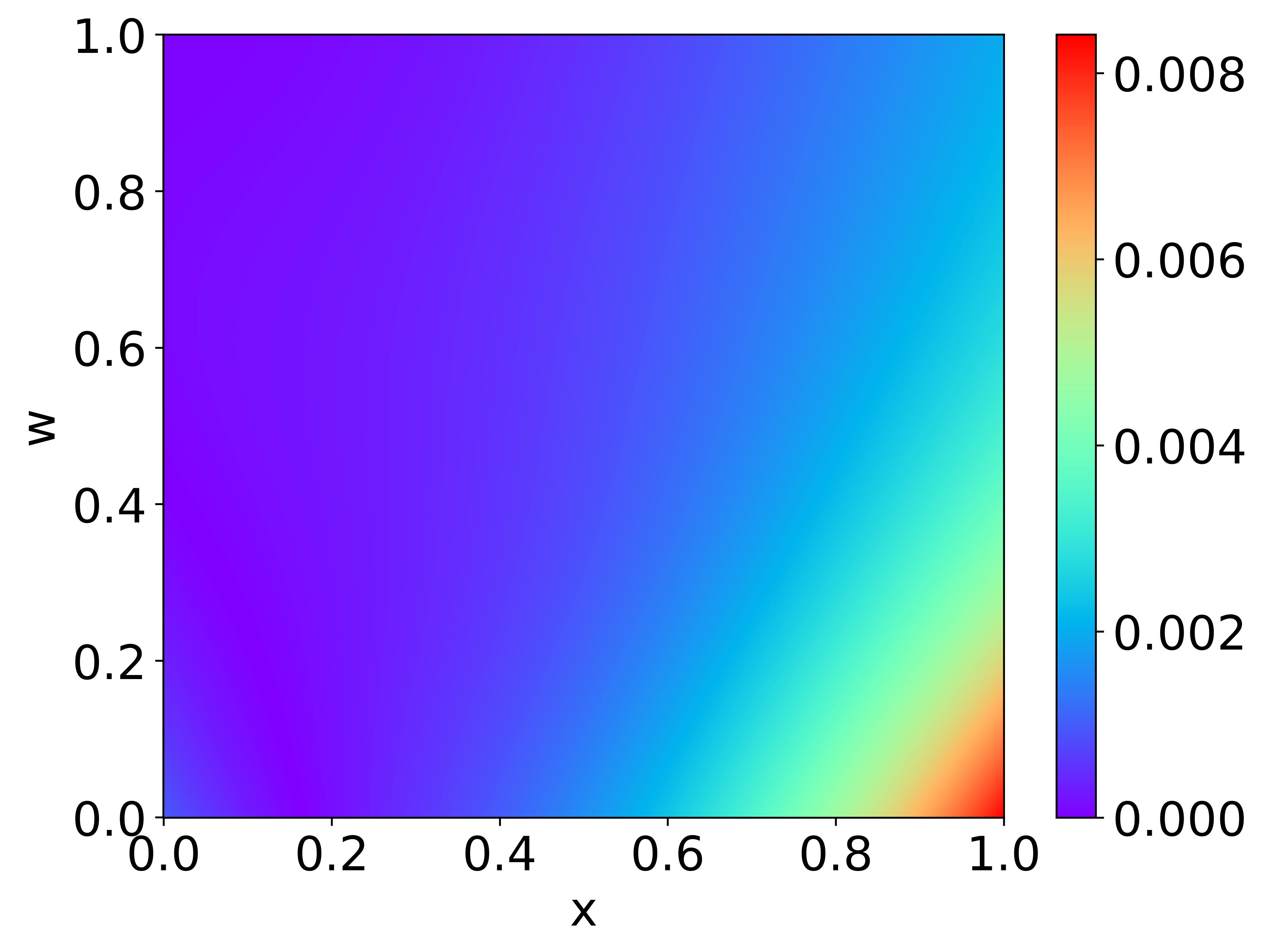}}
\\
\subfigure[absolute error of DGM-PIA value function  ]{\includegraphics[width=0.45\textwidth]{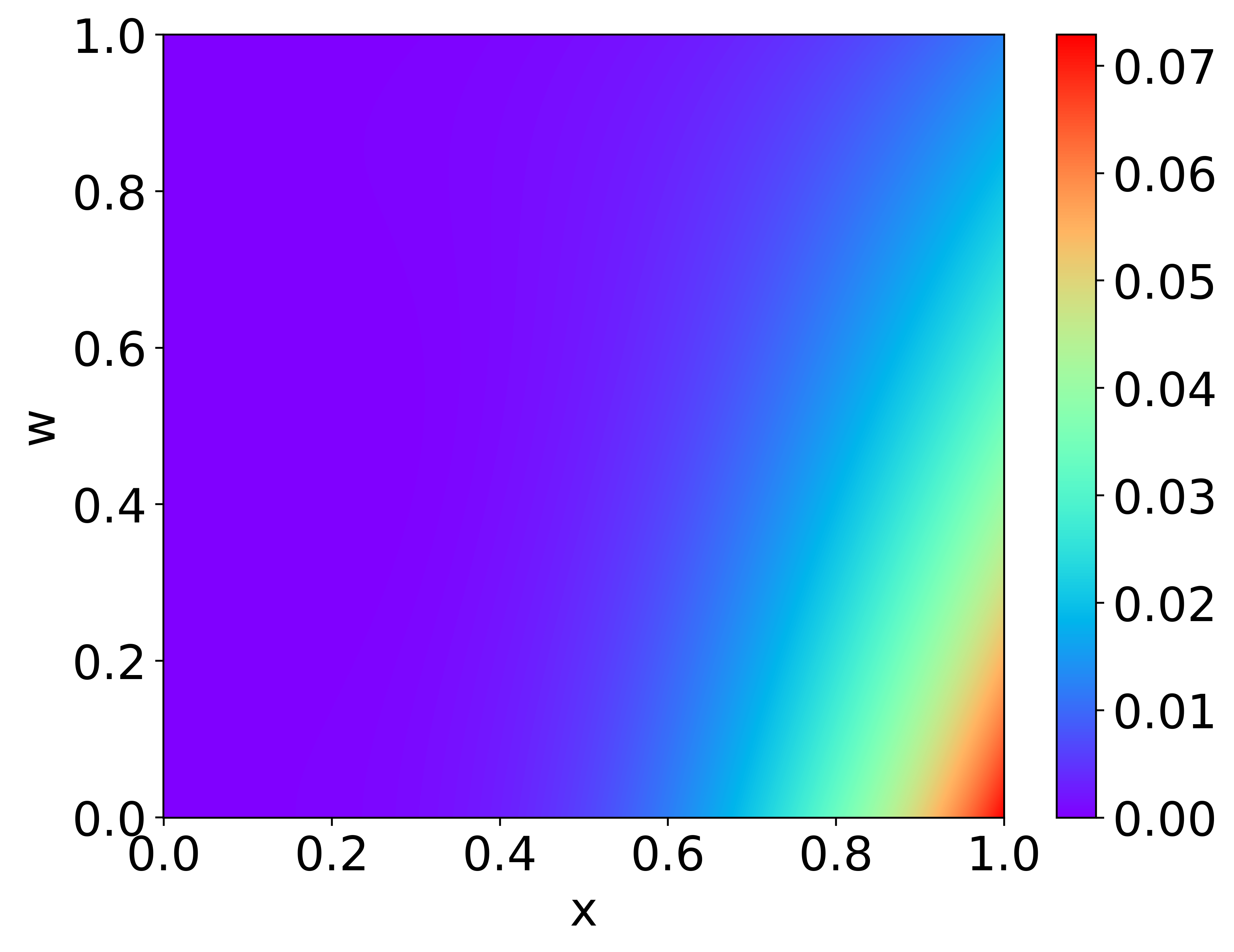}} 
\subfigure[absolute error of DGM-PIA optimal control ]{\includegraphics[width=0.45\textwidth]{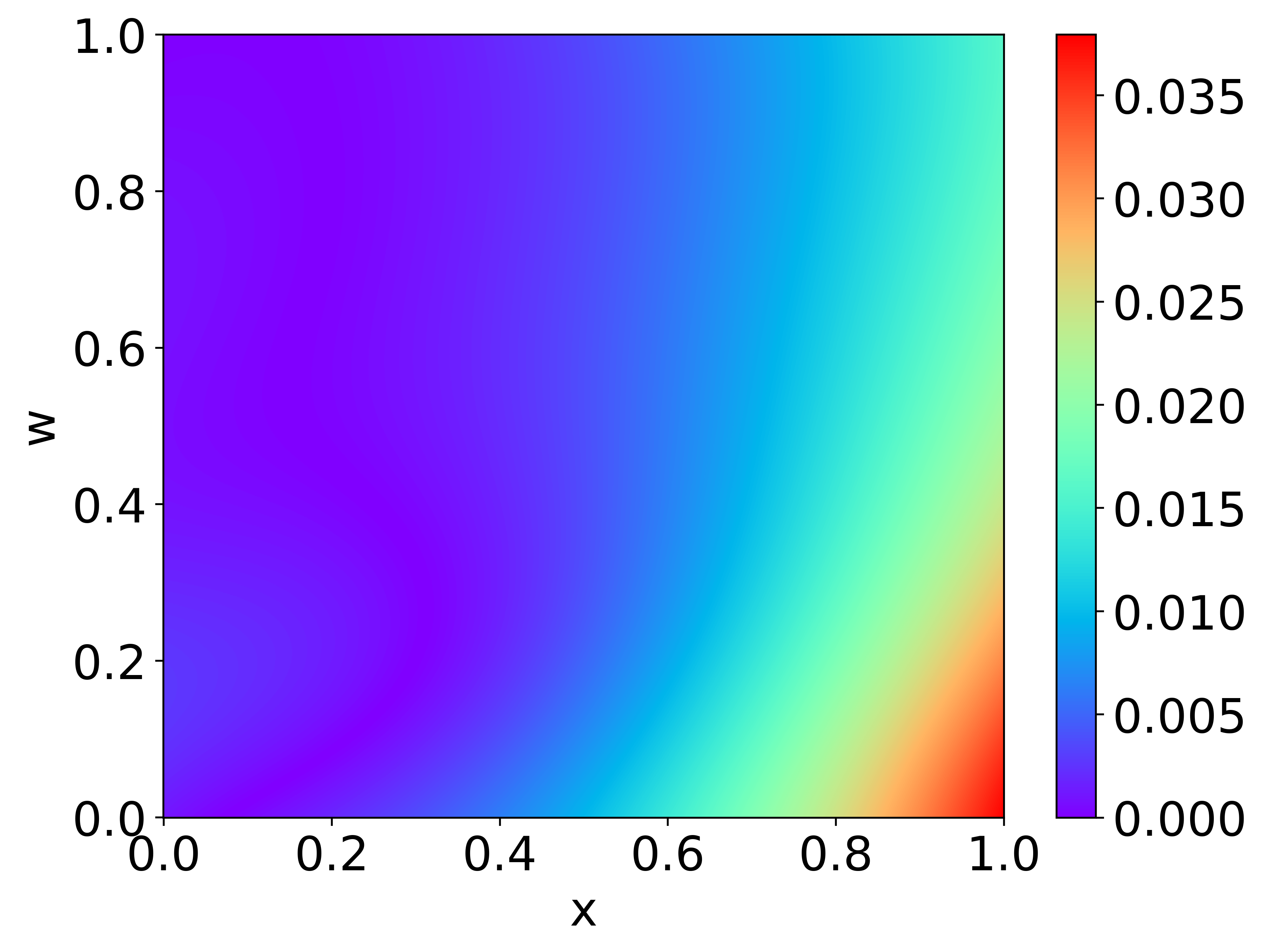}}
\caption{ Heatmaps of the absolute errors in the value function $V(0,x,w)$ and optimal control $a(0,x,w)$ at $t=0$ relative to the ground truth for the case study of Section \ref{sec:mul}. \emph{Top row:}  DeepPAAC algorithm. \emph{Bottom row:}  DGM-PIA  from \cite{al2019applications,al2022extensions}. \label{fig:Mul-explicit}}
\end{figure}

For the NN training, we use $M=2000$ training points per epoch and $B=10$ SGD steps. The learning rate schedule is of polynomial type with power $0.8$, $\alpha_0 =10^{-3}$ decaying to $\alpha_{10^5} = 10^{-5}$. With the stopping rule $\mathrm{Tol_{int}}=\mathrm{Tol_{ctrl}}=10^{-3}$, convergence is achieved in $51570$ training iterations for our DeepPAAC algorithm, whereas the DGM algorithm converges only after $94380$ training iterations, see Figure \ref{fig:Mul-explicit}. Moreover, even after convergence the resulting absolute error compared to the exact ground truth solution are much worse: e.g.~for $(x,w)$ in the neighborhood of $(1,0)$ (bottom right corner in Figure \ref{fig:Mul-explicit}) the Principal's value function of DGM-PIA has error of $>0.06$ compared to just $0.0025$ for DeepPAAC.

\subsection{Multidimensional Problem with Multidimensional Controls}\label{sec:CPT}

As a final, and most numerically challenging example, we consider a two-dimensional state and two-dimensional controls.  In Section \ref{subsec: multi control-one dim spatial}, the examples had multidimensional controls, but there was just one spatial dimension. In this section, we numerically solve the example from Section 5.3 in Cvitani{\'c}, Possama{\"\i} and Touzi \cite{cvitanic2018dynamic}. This model features an Agent who controls both the drift and the volatility of $\{X_t\}$:
we have  $b(t,a,x)=a$,  while the volatility term also depends on $a$: $\sigma(t,a,x)=a$\footnote{ When there is volatility control in the agent's problem, the well-posedness of the PA  problem is much more complicated. Cvitani{\'c} et al.~\cite{cvitanic2018dynamic} utilize second order BSDEs to establish the well-posedness of this specific model. Although this example falls beyond the framework of Section \ref{sec:model},  we nevertheless use it as one of our numerical tests, since our work focuses more on the numerical implementation.}. There are no running costs for either player
$$U^A(t,\alpha,\beta,a,c)=U^P(t,\alpha,\beta,D)=0.$$
For the lump sum payment $\xi$, the  Agent has a risk-averse utility
$\Phi^A(T,\xi,m)=\Phi^A(\xi)$, and so does the Principal, $\Phi^P(T,\xi,X,m)=\Phi^P(X-\xi)$ for some concave functions $\Phi^A$ and $\Phi^P$ to be specified later.

The main difference in this example is that the Agent's problem has no cost on $a_t$, therefore the problem becomes \emph{first-best}, i.e., the Principal effectively controls the Agent's effort and the Principal's problem features two controls.  The Agent's wealth process $\{W_t\}$ evolves as
$$dW_t=a_tZ_t \, dB_t^a,  \qquad W_T=\Phi^A(\xi).$$

After a certain transformation, \cite{cvitanic2018dynamic} views $Z$ as ${\partial V^A\over \partial x}$ and $-a$  as ${\partial^2 V^A\over \partial x^2}$ (Recall $V^A$ is Agent's value function). The HJB equation for the Principal's problem is therefore
\begin{eqnarray}
\begin{cases}
\label{eq:CPT example}
V^P_t+\sup\limits_{Z,a} \left\{a V^P_x+{1\over2}a^2(V^P_{xx}+Z^2V^P_{ww})+a^2ZV^P_{xw}\right\} =0,\\
V^P(T,x,w)=\Phi^P \left(x-(\Phi^A)^{-1}(w) \right). 
\end{cases}
\end{eqnarray}

The Hamiltonian is quadratic in each of $a$ and $Z$ however note the high-order term $\frac{1}{2} a^2 Z^2 V^P_{ww}$. This yields the 2-dimensional  loss criterion (taking the gradient in $a,Z$ respectively):
\begin{align*}
L_{ctrl}=\Big[V^P_x +a(V^P_{xx}+Z^2V^P_{ww})+ 2aZV^P_{xw}, \; a^2zV^P_{ww}+a^2V^P_{xw}\Big].
\end{align*}

In our numerical example, we choose exponential terminal utilities $\Phi^P(x) = \exp(-\gamma_P x)$, $\Phi^A(x) = \exp(-\gamma_A x)$, which lead to the terminal condition 
$$V^P(T,x,w)=-e^{-\gamma_P x}{1\over (-w)^{\gamma_P\over \gamma_A}}. $$

\begin{figure}[!htbp]
\centering
\subfigure[$V^P(0,x,w)$ at $t=0$]{\includegraphics[width=0.32\textwidth]{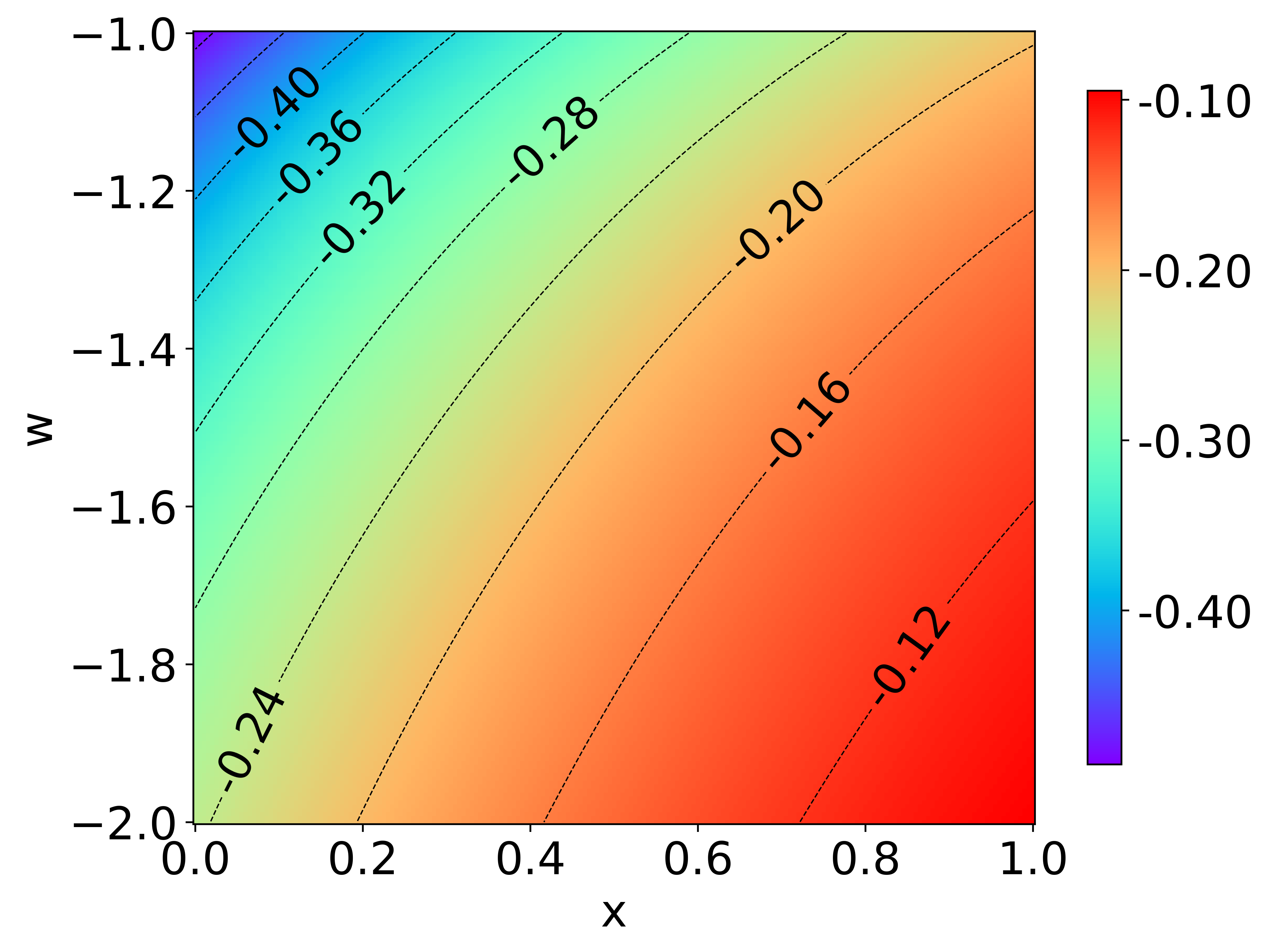}} 
\subfigure[$Z(0, x,w)$ at $t=0$]{\includegraphics[width=0.32\textwidth]{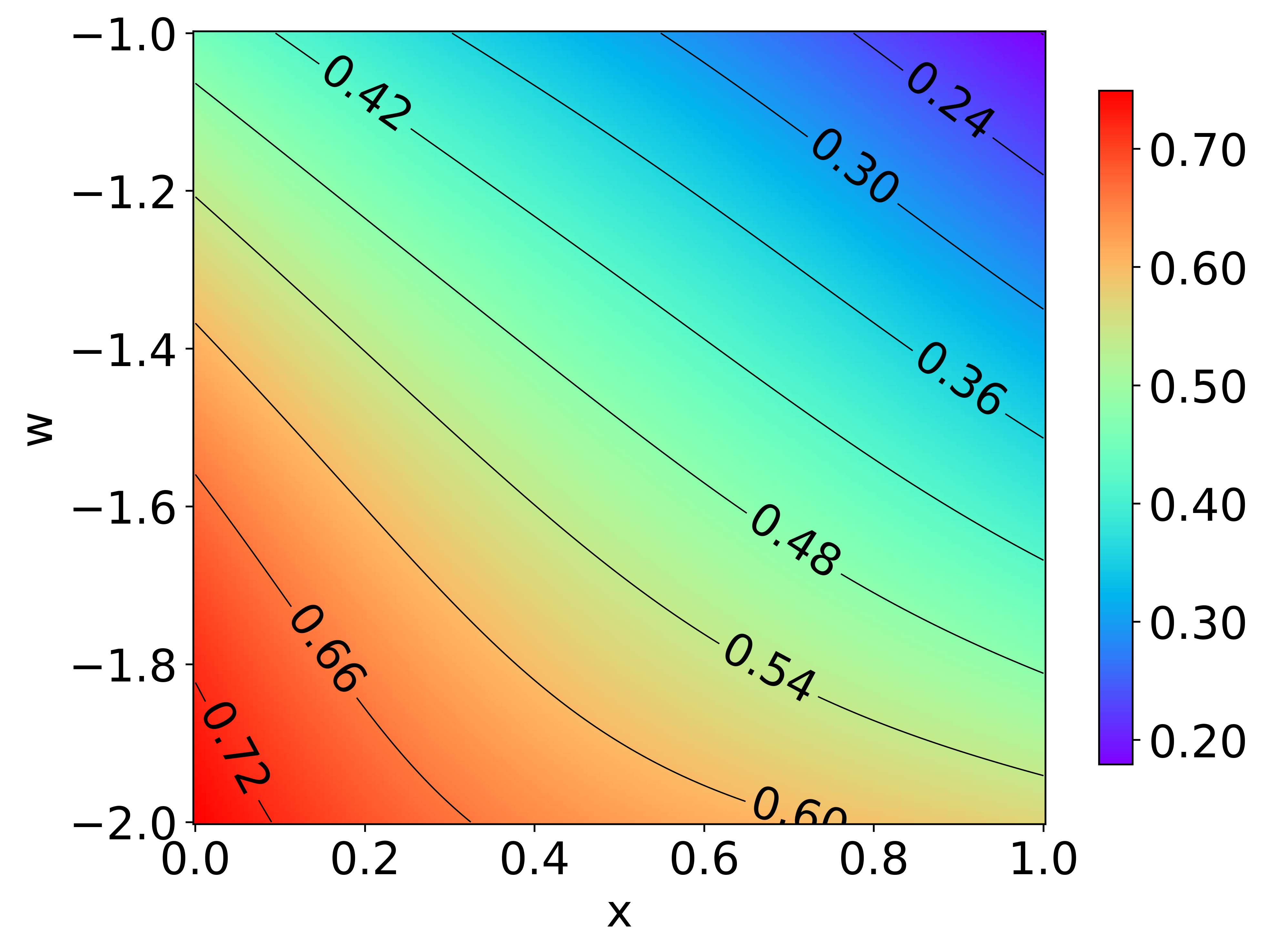}}
\subfigure[$a(0, x,w)$ at $t=0$]{\includegraphics[width=0.32\textwidth]{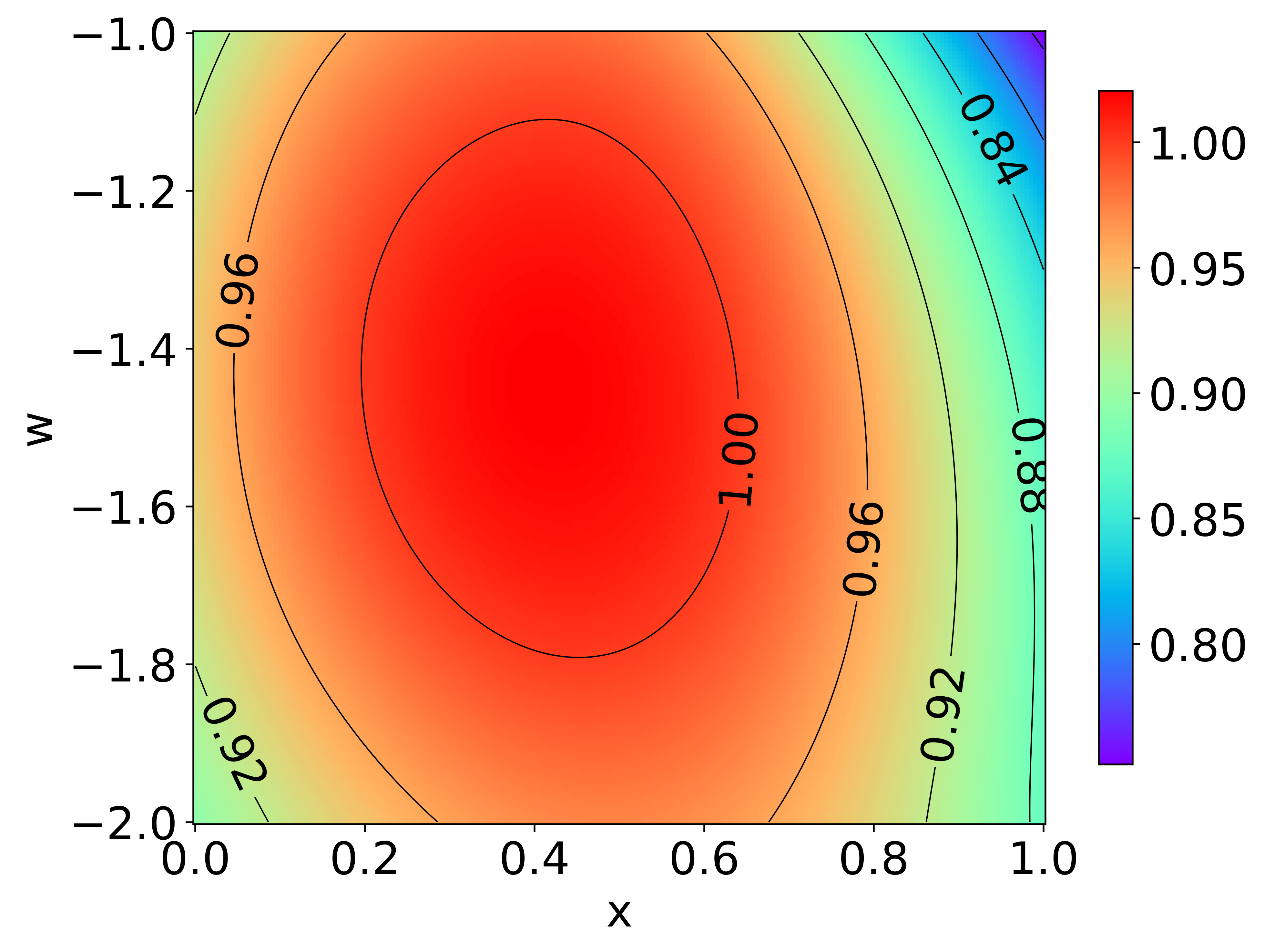}} \\
\subfigure[$V^P(0,x,w)$ at $t=0$]{\includegraphics[width=0.32\textwidth]{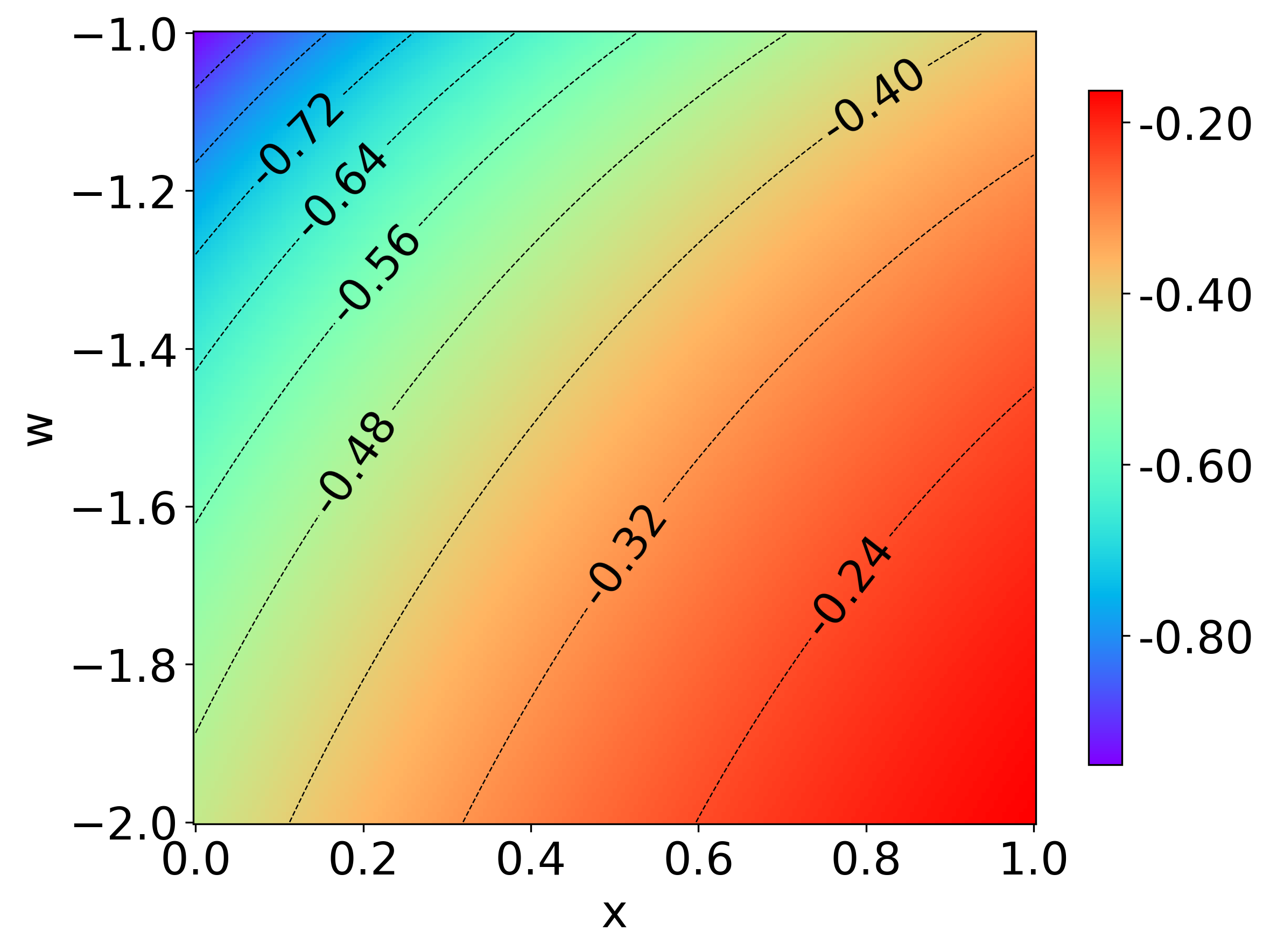}} 
\subfigure[$Z(0, x,w)$ at $t=0$]{\includegraphics[width=0.32\textwidth]{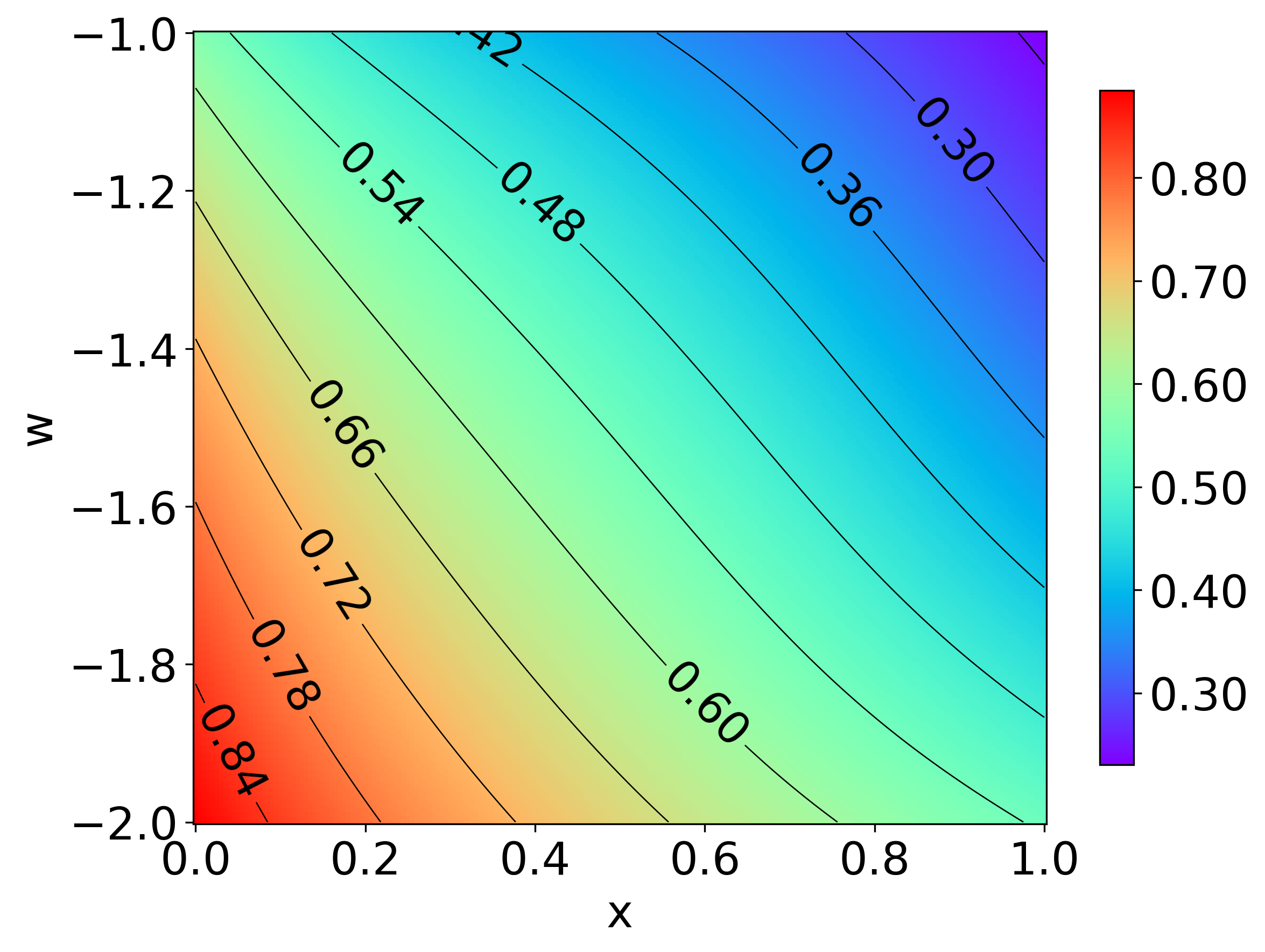}}
\subfigure[$a(0, x,w)$ at $t=0$]{\includegraphics[width=0.32\textwidth]{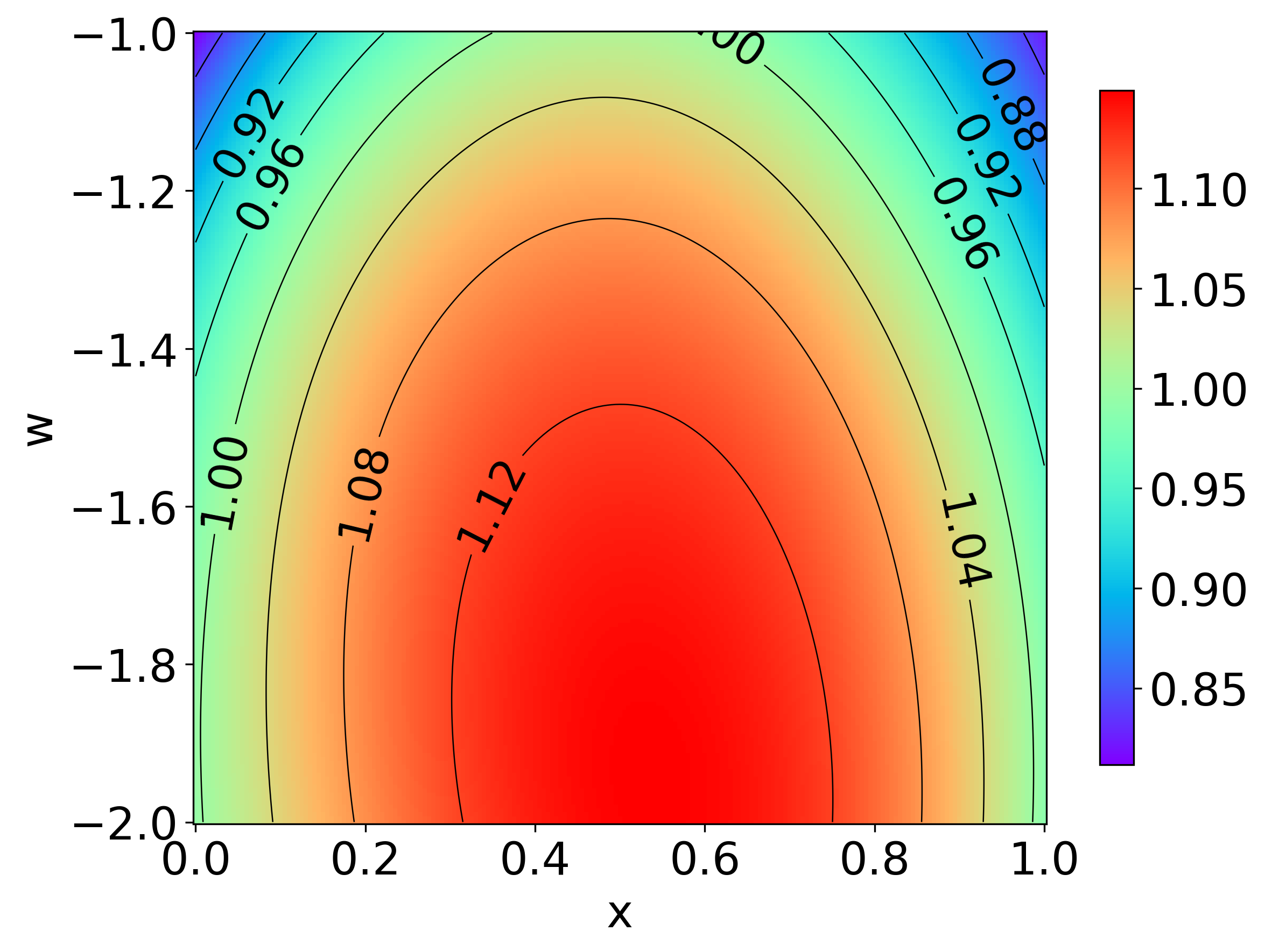}}
\caption{Value function and optimal controls of the example of Section \ref{sec:CPT}. Panels (a)-(c):  $\gamma_A = \gamma_P = 1$ (after convergence at $n=85230$ steps). 
Panels (d)-(f): $\gamma_A = 2, \gamma_P = 1$ (after convergence at $n=92580$ steps). \label{fig:2D-example-vo} }
\end{figure}

The equation \eqref{eq:CPT example} does not have an explicit solution. 
Figure \ref{fig:2D-example-vo} shows the solutions for the cases $\gamma_P = \gamma_A = 1$ as well as $\gamma_P=1, \gamma_A=2$. To evaluate the first and second derivatives of $V^P$ we use the built-in autodiff functions of the neural network library. 
To solve this problem, we take  the convergence thresholds $\mathrm{Tol_{int}}=\mathrm{Tol_{ctrl}}=10^{-3}$, size of the training set is $M=2000$, with  $B=30$ steps per epoch, and polynomial learning rate decay from $10^{-3}$ to $10^{-4}$ with power $0.8$. Numerical convergence is achieved in $85230/92580$ training iterations. 

In Figure \ref{fig:2D-example-vo} (a)-(d), it is seen that the Principal's value function $V^P$ increases from upper left to lower right, meaning that it is increasing in $x$ and decreasing in $w$. This is consistent with the theoretical result that the Principal's value function increases when the initial value of output process $\{X_t\}$ increases or the Agent's participant constraint $R$ decreases. Figures \ref{fig:2D-example-vo}(b)-(e) show that $Z={\partial V^A\over \partial x}$ decreases in both $x$ and $w$, which is consistent with the concavity condition for the Agent's value function: ${\partial^2 V^A\over \partial x^2}<0$. The latter also implies that ${\partial^2 V^A\over \partial x \partial w}<0$ for this example. Figure \ref{fig:2D-example-vo}(c) shows that the Agent's optimal effort process $a^*$ (or equivalently, ${\partial^2 V^A\over \partial x^2}$) has a peak in the middle of the state space, maximized at moderate positive $x \simeq 0.5$ and moderate negative $w \simeq -1.5$. This peak shifts to a lower $w$ when we change to $\gamma_P = 2$ (panel f).

\section{Conclusion}
\label{sec:Con}

In this article, we have developed and validated a deep learning algorithm for Principal-Agent problems. We employ an actor-critic setup to solve the corresponding HJB PDE with an implicit Hamiltonian. While so far most PA setups have been very restrictive, through a collection of case studies we demonstrate that the proposed DeepPAAC algorithm is able to handle many extensions, from multi-dimensional state/multi-dimensional controls, to (coupled) control constraints. The new scheme opens the door for future numerical-driven investigations of novel PA settings. In tandem, it also sheds light on aspects of using neural networks for solving HJB equations, an ongoing active area of research across applied mathematics. \\

{\bf Acknowledgement}: The authors would like to thank two anonymous referees whose constructive comments have greatly helped to improve the paper. The authors would also like to thank Chenchen Mou and Yuri Saporito for fruitful discussions.  Mike Ludkovski is partially supported by NSF DMS-2407550 and DMS-2420988. Zimu Zhu is supported by The Hong Kong University of Science and Technology(Guangzhou) Start-up Fund G0101000240 and the Guangzhou-HKUST(GZ) Joint Funding Program (No. 2024A03J0630).

\section{Appendix: Proof of Proposition \ref{prop:HMmodel}}
\proof We prove the result in 2 steps:

{\bf Step 1}: Solving the Agent's problem: for $0\leq t\leq T$, $\Phi^A(x)=-\exp(-\gamma_A x)$, set 
\begin{eqnarray*}
W_t&=\sup_{a}(\Phi^A)^{-1}\left(\EE^a \left[\Phi^A(\xi-{1\over 2}\int_t^T a_s^2ds)\big{|}\,\mathcal{F}_t \right]\right).
\end{eqnarray*}

Using \^{I}to's formula and the comparison principle of BSDE (See e.g., \cite{zhang2017backward}), we have that $W_t$ satisfies the following BSDE:
\begin{align}
\label{eq:Y process for HM model}
W_t=\xi-\int_t^T {\gamma_A-1\over 2}Z_s^2ds-\int_t^T Z_s dB_s\quad 0\leq t\leq T,
\end{align}
\noindent and the Agent's optimal effort process is $a^*_t=Z_t, \forall 0\leq t\leq T$. When the Agent's participation constraint is $R_t$ at time $t$, we have the initial condition $W_t=-{1\over \gamma_A}\log(-R_t)$.

{\bf Step 2:}  Solving the Principal's problem: we first rewrite (\ref{eq:Y process for HM model}) into the forward form:
\begin{align*}
\xi&=W_t+\int_t^T {\gamma_A-1\over 2}Z_s^2 ds+\int_t^T Z_s dB_s\\
&=W_t+\int_t^T {\gamma_A+1\over 2}Z_s^2 ds+\int_t^T Z_s dB_s^Z.
\end{align*}
Therefore, the Principal's maximization problem is (since $a_t^*=Z_t$)
\begin{eqnarray*}
&&\sup_{Z} \EE^Z \left[\Phi^P \left(\int_t^T Z_sds+\int_t^T dB_s^Z-Y_t-\int_t^T {\gamma_A+1\over 2}Z_s^2ds-\int_t^T Z_s dB_s^Z \right)\big{|} \, \mathcal{F}_t  \right]\\
 =&& \sup_{Z}\EE^Z \left[\Phi^P \left(-W_t+\int_t^T (Z_s-{1+\gamma_A\over 2}Z_s^2) ds+\int_t^T (1-Z_s)dB_s^Z \right) \big{|} \, \mathcal{F}_t \right].
\end{eqnarray*}
We re-interpret above as the problem of maximizing 
$\sup\limits_{Z}  \EE^Z \left[\Phi^P(\tilde{X}_T) \right]$
where the auxiliary state process follows
\begin{align*}
d\tilde{X}_s=\big(Z_s-{1+\gamma_A\over 2}Z_s^2 \big)ds+(1-Z_s)dB_s^Z,\quad t\leq s\leq T, \qquad \tilde{X}_t=x.
\end{align*}
 When Agent's participant constraint is $R_t$ at time $t$, the initial value is $x={1\over \gamma_A}\log(-R_t)$. This yields  the HJB equation for the Principal:
\begin{equation}
\label{eq:appendix-HMHJB}
\left\{
\begin{aligned}
&V^P_t+\sup_Z \left\{ V^P_x(Z-{1+\gamma_A\over 2}Z^2)+{1\over 2}V^P_{xx}(1-Z)^2)\right\}=0;\\
&V^P(T,x)=\Phi^P(x).\\
\end{aligned}
\right.
\end{equation}
When $\Phi^P(x)=-\exp(-\gamma_P x)$, we make the ansatz for the value function as $V^P(t,x)=\Phi^P\big(x+c(T-t)\big)$ for a constant $c$ to be determined. Plugging this ansatz  into (\ref{eq:appendix-HMHJB}), we get
\begin{align*}
V^P(t,x)=\Phi^P \left(x+\Bigl[{1\over 2}{(1+\gamma_P)^2\over 1+\gamma_A+\gamma_P}-{\gamma_P\over 2} \Bigr](T-t) \right)
\end{align*}
which implies \eqref{eq:a-HM} and \eqref{eq:xi-HM}.
\qed

\bibliographystyle{siam}
\bibliography{bibiography}

\end{document}